\documentclass[preprint]{imsart} 

\RequirePackage[OT1]{fontenc}
\RequirePackage{amsthm,amsmath}
\RequirePackage[numbers]{natbib}
\RequirePackage[colorlinks,citecolor=blue,urlcolor=blue]{hyperref}
\usepackage{amsmath}
\usepackage[english]{babel} 
\usepackage{graphicx} 
\usepackage{epstopdf}
\usepackage{float}
\usepackage{fixmath}
\usepackage{bm}
\usepackage{inputenc,t1enc}
\usepackage{bbm}
\usepackage{amsfonts}
\usepackage{latexsym}
\usepackage{mathrsfs}
\usepackage{mathtools,xparse}
\usepackage{enumitem}  
\usepackage{subcaption}

\firstpage{0}
\lastpage{0}   

\DeclarePairedDelimiter{\norm}{\lVert}{\rVert}

\startlocaldefs
\numberwithin{equation}{section}
\theoremstyle{plain}
\newcommand{\nc}{\newcommand}
\nc{\nt}{\newtheorem}
\nt{defn}{Definition}
\nt{lem}{Lemma}
\nt{pr}{Proposition}
\nt{theorem}{Theorem}
\nt{cor}{Corollary}
\nt{remark}{Remark}
\nt{ex}{Example}
\nt{ass}{Assumption}
\nt{step}{Step}
\nt{case}{Case}
\nt{subcase}{Subcase}
\nt{note}{Note}
\nt{assum}{Assumption}
\nc{\bd}{\begin{defn}} \nc{\ed}{\end{defn}}
\nc{\blem}{\begin{lem}} \nc{\elem}{\end{lem}}
\nc{\bpr}{\begin{pr}} \nc{\epr}{\end{pr}}
\nc{\bth}{\begin{theorem}} \nc{\eth}{\end{theorem}}
\nc{\bcor}{\begin{cor}} \nc{\ecor}{\end{cor}}
\nc{\bex}{\begin{ex}}  \nc{\eex}{\end{ex}}
\nc{\bass}{\begin{ass}}  \nc{\eass}{\end{ass}}
\nc{\bstep}{\begin{step}}  \nc{\estep}{\end{step}}
\nc{\bcase}{\begin{case}}  \nc{\ecase}{\end{case}}
\nc{\bsubcase}{\begin{subcase}}  \nc{\esubcase}{\end{subcase}}
\nc{\bnote}{\begin{note}}  \nc{\enote}{\end{note}}
\nc{\bassum}{\begin{assum}}  \nc{\eassum}{\end{assum}}
\nc{\prf}{{\bf Proof.} }
\nc{\eop}{\hfill $\Box$ \\ \\}
\nc{\argmin}{\mathrm{argmin}}
\nc{\argmax}{\mathrm{argmax}}
\nc{\sgn}{\mathrm{sgn}}
\nc{\Var}{\mathrm{Var}}
\nc{\Cov}{\mathrm{Cov}}
\nc{\bak}{\!\!\!\!\!}
\nc{\IBD}{\mathrm{IBD}}
\nc{\supp}{\mathrm{supp}}
\nc{\dom}{\mathrm{dom}}
\nc{\R}{{\mathbb R}}
\nc{\peq}{\preceq}
\nc{\wt}{\widetilde}
\nc{\Mult}{\mathrm{Mult}}

\nc{\Prob}[1]{\mathbb{P}_{#1}}

\newcommand{\tcb}{\textcolor{black}}

\endlocaldefs

\begin{document}

\begin{frontmatter}

\title{The asymptotic distribution of the isotonic regression estimator over a general countable pre-ordered set}
\runtitle{Isotonic regression over a countable pre-ordered set}

\begin{aug}
\author{\fnms{Dragi} \snm{Anevski}\ead[label=e1]{dragi@maths.lth.se}}
\and
\author{\fnms{Vladimir} \snm{Pastukhov}\corref{craut1}\ead[label=e2]{pastuhov@maths.lth.se}}

\address{Center for Mathematical Sciences\\
Lund University\\
\printead{e1,e2}}
\runauthor{D. Anevski and V. Pastukhov}
\end{aug}

\begin{abstract}
We study the isotonic regression estimator over a general countable pre-ordered set. We obtain the limiting distribution of the estimator and study its properties. It is proved that, under some general assumptions, the limiting distribution of the isotonized estimator is given by the concatenation of the separate isotonic regressions of the certain subvectors of an unrestrecred estimator's asymptotic distribution. Also, we show that the isotonization preserves the rate of convergence of the underlying estimator. We apply these results to the problems of estimation of a bimonotone regression function and estimation of a bimonotone probability mass function.
\end{abstract}




\tableofcontents


\end{frontmatter}

\section{Introduction}

In this paper we study estimators for the problem of estimating real valued functions that are  defined on a countable set and that are monotone with respect to a pre-order defined on that set. In the situation when there exists an underlying, empirical, estimator of the function, which is not necessarily monotone, but for which one has (process) limit distribution results, we are able to provide limit distribution results  for the estimators. Our results can be applied to the special cases of probability mass function (pmf) estimation and regression function estimation.  In the case of estimating a bimonotone pmf, i.e. a pmf which is monotone with respect to the usual matrix pre-order on ${\mathbb Z}^2_+$, we state the limit distribution of the order restricted maximum likelihood estimator (mle), thereby generalising previously obtained results by \cite{jankowski2009estimation}, who treated the one-dimensional case, i.e. the mle of a monotone pmf on ${\mathbb Z}_+$. In fact we are able to state limit distribution results for the mle of a monotone pmf on ${\mathbb Z}_+^d$, for arbitrary $d>1$, cf. Corollary \ref{dmonot} below. In the case of estimating a bimonotone regression function, i.e. a function defined on ${\mathbb Z}_+^2$ that is monotone with respect to the matrix pre-order on ${\mathbb Z}_+^2$, we state the limit distribution of the isotonic regression estimator, again generalising previously known results for the isotonic regression on ${\mathbb Z}_+$, cf. \cite{jankowski2009estimation}. In this setting we would also like to mention \cite{dumbgen}, that studied algorithms resulting from the minimisation of a smooth criterion function under bimonotonicity constraints. In the regression setting we are able to derive the limit distributions for the isotonic regression of functions that are monotone with respect to the matrix pre-order on ${\mathbb Z}_+^d$, for arbitrary $d>1$, cf. Corollary \ref{dmonotreg}.

We would like to emphasize that the general approach taken in this paper allows for other pre-orders than the usual matrix order on ${\mathbb Z}^d$. Furthermore, our approach allows for also other starting basic empirical estimators; one could e.g. consider non-i.i.d. data settings, treating e.g. stationary (spatially homogenous) dependent data. In fact one can consider our estimator as the final step in a, at least, two-step, approach where in the first, or next-to-last, step, one provides the "empirical" estimator $\hat{\bm{g}}_{n}$ of the estimand $\mathring{\bm{g}}$, for which it is necessary to have established (process) limit distribution result of the form 
\begin{eqnarray}
      n^{1/2}(\hat{\bm{g}}_{n} - \mathring{\bm{g}}) \stackrel{d}{\to} \bm{\lambda} \label{eq:assumption-basic}
\end{eqnarray}
on the appropriate space, e.g. $l^2$, see Assumptions \ref{astheta} and \ref{asthetainf} below. Note that we have simplified Assumptions  \ref{astheta} and \ref{asthetainf}  slightly in $(\ref{eq:assumption-basic})$ for illustrative purposes; the rate $n^{1/2}$ in $(\ref{eq:assumption-basic})$ is allowed to differ, even between the components in the vector $\hat{\bm{g}}_{n}$. Given the assumption  $(\ref{eq:assumption-basic})$ we then establish a limit distribution result for the proposed estimator $\hat{\bm{g}}^{*}_{n}$ of the form
\begin{eqnarray}
  n^{1/2}(\hat{\bm{g}}_{n}^{*} - \mathring{\bm{g}}) \stackrel{d}{\to} \varphi(\bm{\lambda}),
\end{eqnarray}
in Theorems \ref{thmasymfin} and \ref{thmasyminf}, where $ \varphi$ is a certain isotonic regression operator defined in the sequel.

The general approach in this paper is somewhat reminiscent to the approach taken in \cite{anevhos}, in which one considered a two-step general procedure for isotonization, allowing  e.g. different types of dependence structures on the data. The difference to our paper is that we treat the, arguably, more complex notion of monotonicity with respect to pre-orders on $d$-dimensional spaces, whereas \cite{anevhos} only treated monotonicity in the one dimensional setting, and, furthermore,  that we treat only functions with discrete or countable support, such as pmfs, whereas \cite{anevhos} treated functions with continuous support, such as pdfs.

This work is mainly motivated by the results obtained in \cite{baljan, jankowski2009estimation}. In \cite{jankowski2009estimation} the problem of estimation of a discrete monotone distribution was studied in detail. It was shown that the limiting distribution of the constrained mle of a pmf is a concatenation of the isotonic regressions of Gaussian vectors over the periods of constancy of the true pmf $\bm{p}$, cf. Theorem 3.8 in \cite{jankowski2009estimation}. In the derivation of the limiting distribution in \cite{jankowski2009estimation} the authors used the strong consistency of the empirical estimator of $\bm{p}$ as well as the fact that the constrained mle is given by  the least concave majorant (lcm) of the empirical cumulative distribution function (ecdf). 

The problem of maximum likelihood estimation of a unimodal pmf was studied in \cite{baljan}. That problem is different from the one being considered here, since \cite{baljan} treats only pmfs on $\mathbb{Z}$, whereas we are able to treat multivariate problems with our approach.  

In our work we do not require strong consistency of a basic estimator $\hat{\bm{g}}_{n}$, and we consider general pre-order constraints, resulting in an expression for the isotonic regression that is more complicated than the lcm of the ecdf, c.f  Assumptions \ref{astheta} and \ref{asthetainf}. Also it turns out that the limiting distribution of the isotonized estimator $\hat{\bm{g}}_{n}^{*}$ can be split deeper than to the level sets of $\mathring{\bm{g}}$, which are the analogues of the periods of constancy of $\mathring{\bm{g}}$ in the univariate case.

For a general introduction to the subject of constrained inference we refer to the monographs: Barlow R. E. et al. \cite{barlowstatistical},  Robertson T. et al. \cite{robertsonorder}, Silvapulle M. J. \cite{silvapsen} and Groeneboom P. et al. \cite{groeneboom}. In these monographs the problem of an isotonic regression has been considered in different settings, and, in particular, basic questions such as existence and uniqueness of the estimators have been addressed. In Lemmas \ref{propisot} and  \ref{propisotinf} below we list those properties which will be used in the proofs of our results.

There are (at least) two reasons why one may use an isotonic regression. The first, and most important, reason is that one may desire an order restricted function or parameter. Thus, viewing the inference problem as an order restricted problem, one needs to come up with an algorithm that produces an order restricted estimator. A second reason is the following error reduction property: if $\hat{\bm{g}}$ is any vector and if 
$\hat{\bm{g}}^{*}$ is the isotonic regression of $\hat{\bm{g}}$ with weights $\bm{w} = (w_{1}, w_{2}, \dots)$ then 
\begin{eqnarray}\label{erredpr}
\sum_{i} \Phi[\hat{g}_{i}^{*} - h_{i}]w_{i} \leq \sum_{i} \Phi[\hat{g}_{i} - h_{i}]w_{i},
\end{eqnarray}
for any convex function $\Phi$ on $(-\infty, \infty)$ and any isotonic vector $\bm{h}$, cf. Theorem 1.6.1 in \cite{robertsonorder}. Therefore, for any estimator $\hat{\bm{g}}_{n}$ such that $\hat{\bm{g}}_{n} \stackrel{p}{\to} \mathring{\bm{g}}$, where $\mathring{\bm{g}}$ is isotonic, one has
\begin{eqnarray*}
||\hat{\bm{g}}_{n}^{*} - \mathring{\bm{g}}||_{2} \leq ||\hat{\bm{g}}_{n} - \mathring{\bm{g}}||_{2},
\end{eqnarray*}
for any finite $n$, which also immediately implies that $\hat{\bm{g}}_{n}^{*}$ is consistent.

In the paper \cite{jankowski2009estimation} the authors also studied the finite sample behaviour of the isotonized estimator of one-dimensional empirical pmf. It was shown there that the isotonized estimator performs better than empirical estimator in $l_{1}$ and $l_{2}$ metrics and in Hellinger distance. The gain from isotonization depends on the structure of the underlying pmf $\bm{p}$: the more constant regions the true pmf $\bm{p}$ has the better performance of the isotonic regression.

To demonstrate the performance of the isotonic regression in two dimensions, let us consider the following bimonotone pmf $\bm{p}$ on a finite grid $\mathcal{X} := \{ \bm{x}= (i_{1} , i_{2})^{T}: i_{1} = 1,\dots, 5,  i_{2} = 1,\dots, 5\}$
\begin{eqnarray*}
\bm{p} = q_{1}U(1) + \dots + q_{5}U(5),
\end{eqnarray*}
where $\bm{q} = (0.1, 0.2, 0.3, 0.2, 0.2)$ and $U(r)$ is the uniform distribution on $\mathcal{X}_r=\{ \bm{x}= (i_{1} , i_{2})^{T}: i_{1} = 1,\dots, r,  i_{2} = 1,\dots, r\}$. 

\begin{figure}[] 
\includegraphics[scale=0.55]{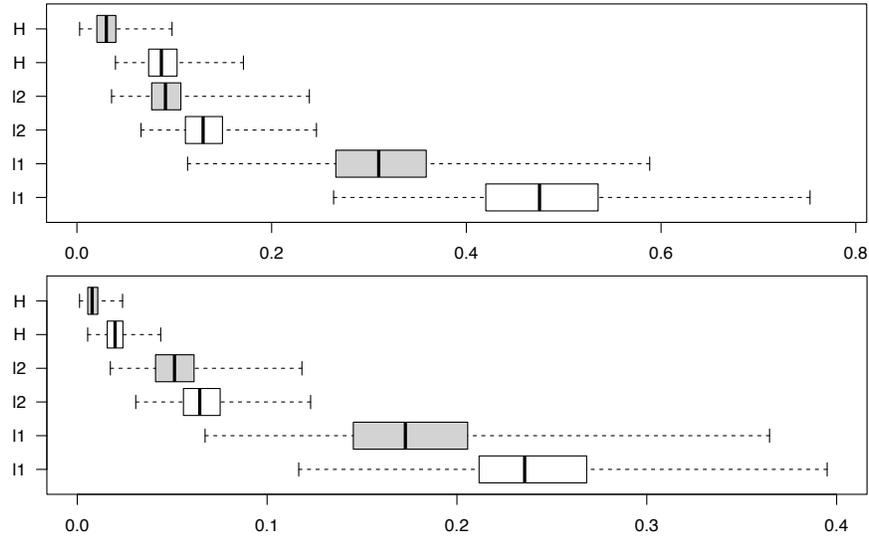}
\caption{Boxplots for $l_1, l_2$ and Hellinger distances for the empirical estimator (white) and isotonic regression (grey) of $\bm{p}$  for the sample sizes $n=50$ (top) and for $n=300$ (bottom). (Note, that the scales for $n=50$ and $n=300$ are different.)}\label{decrconstpmf}
\end{figure}
Figure \ref{decrconstpmf} shows the boxplots of $1000$ samples for $l_1, l_2$ and Hellinger distances for the empirical estimator and for the isotonic regression  of $\bm{p}$, respectively, for the sample sizes $n=50$ at the top and for $n=300$ at the bottom. One can see that the isotonizatized estimator performs significantly better than the empirical estimator.

The asymptotic behaviour of the regression estimates over a continuous setup under monotonic restriction was first studied in \cite{brunk70, wright81}, where it was shown that the difference of the regression function and its estimate multiplied by $n^{1/3}$, at a point with a positive slope, has a nondegenerate limiting distribution.  The problem of estimating a monotone pdf was studied, for example, in \cite{anevhos, carolandyks, grenander, rao}. In \cite{anevhos} the authors studied a general scheme for order constrained inference in a continuous setup. In the discrete case some recent results are \cite{balabdurkol, baljan, sengun, jankowski2009estimation}. In  \cite{sengun} the authors studied  risk bounds for isotonic regression.

The remainder of this paper is organised as follows. In Section \ref{sec:definitions} we introduce some notations and define the estimators. In Section \ref{sec:finitecase} we consider the finite dimensional case. Theorem \ref{thmasymfin} gives the asymptotic distribution of the isotonized estimator. Next, in Section \ref{sec:infinitecase} we consider the infinite dimensional case, which is quite different from the finite one. Theorem \ref{thmasyminf} describes the asymptotic behaviour of the isotonized estimator for the infinite dimensional case. In Section \ref{sec:applications} we first discuss the application of the obtained results to the problems of estimation of a bimonotone regression function and of a bimonotone probability mass function, respectively, \tcb{and then the corresponding limit distribution result for $d$-monotone functions, for an arbitrary $d>1$. The limit distributions are stated in Theorems \ref{bimonotreg}, \ref{bimonot} and Corollaries \ref{dmonotreg} and \ref{dmonot}.}
In Section \ref{sec:discussion} we make some final comments about our results and relations to similar problems. We have gathered proofs of some intermediate results that are stated in the main body of the paper in an Appendix.

\section{The inference problem and notations}\label{sec:definitions}

In order to introduce the inference problem in detail, we start by introducing some notations. Let $\mathcal{X}$ be a countable set $\{ x_{1}, x_{2}, \dots \}$ with $s = |\mathcal{X}| \leq \infty$, with a pre-order $\preceq$ defined on it. We begin with the definitions of the order relations on an arbitrary set $\mathcal{X}$ and of an isotonic regression over it, cf.  also \cite{barlowstatistical, robertsonorder, silvapsen}.
\bd
A binary relation $\preceq$ on $\mathcal{X}$ is a simple order if 
\begin{enumerate}[label=(\roman*)]
\item it is reflexive, i.e. $x \preceq x$ for $x \in \mathcal{X}$;
\item it is transitive, i.e. $x_{1}, x_{2}, x_{3} \in \mathcal{X}$, $x_{1} \preceq x_{2}$ and $x_{2} \preceq x_{3}$ imply $x_{1} \preceq x_{3}$;
\item it is antisymmetric, i.e. $x_{1}, x_{2} \in \mathcal{X}$, $x_{1} \preceq x_{2}$ and $x_{2} \preceq x_{1}$ imply $x_{1} = x_{2}$;
\item it is total, i.e. every two elements of $\mathcal{X}$ are comparable, i.e. $x_{1}, x_{2} \in \mathcal{X}$ implies that  either $x_{1} \preceq x_{2}$ or $x_{2} \preceq x_{1}$.
\end{enumerate}
A binary relation $\preceq$ on $\mathcal{X}$ is a partial order if it is reflexive, transitive and antisymmetric, but there may be noncomparable elements. A pre-order is reflexive and transitive but not necessary antisymmetric and the set $\mathcal{X}$ can have noncomparable elements. Note, that in some literature the pre-order is called as a quasi-order.
\ed

Let us introduce the notation $x_{1} \sim  x_{2}$, if $x_{1}$ and $x_{2}$ are comparable, i.e. if $x_{1} \preceq x_{2}$ or $x_{2} \preceq x_{1}$. 
\bd
A function $f(x): \mathcal{X} \to \mathbb{R}$ is isotonic if $x_{i}, x_{j} \in \mathcal{X}$ and $x_{i} \preceq x_{j}$ imply $f(x_{i}) \leq f(x_{j})$.
\ed

Let $\mathcal{F}^{is} = \mathcal{F}^{is}(\mathcal{X})$ denote the family of real valued bounded functions $f$ on a set $\mathcal{X}$, which are isotonic with respect to the pre-order $\preceq$ on $\mathcal{X}$. In the case when $|\mathcal{X}| = \infty$ we consider the functions from the space $l_{2}^{\bm{w}}$, the Hilbert space of real-valued functions on $\mathcal{X}$, which are square summable with some given non-negative weights $\bm{w} = (w_{1}, w_{2}, \dots )$, i.e. any $g \in l_{2}^{\bm{w}}$ satisfies $\sum_{i=1}^{\infty} g(x_{i})^{2} w_{i} < \infty$. We use the same notation $\mathcal{F}^{is}$ to denote the functions from $l_{2}^{\bm{w}}$ which are isotonic with respect to the pre-order $\preceq$.

\bd
A function $g^{*}: \mathcal{X} \to \mathbb{R}$ is the isotonic regression of a function $g: \mathcal{X} \to \mathbb{R}$ over the pre-ordered set $\mathcal{X}$ with weights $\bm{w} \in \mathbb{R}^{s}_{+}$, with $s\leq \infty$, if 
\begin{eqnarray*}
g^{*} = \underset{f \in \mathcal{F}^{is}}{\argmin} \sum_{x \in \mathcal{X}}(f(x) - g(x))^{2}w_{x},
\end{eqnarray*}
where $w_{x_{i}} = w_{i}$, for $i = 1, \dots, s$.
\ed
Conditions for existence and uniqueness of $g^{*}$ will be stated below.

Similarly one can define an isotonic vector in $\mathbb{R}^{s}$, with $s \leq \infty$, and the isotonic regression of an arbitrary vector in $\mathbb{R}^{s}$. Let us consider a set of indices $\mathcal{I}=\{1, \dots, s\}$, with $s \leq \infty$, with some pre-order $\preceq$ defined on it.
\bd
A vector $\bm{\theta} \in \mathbb{R}^{s}$, with $s\leq \infty$, is isotonic if $i_{1}, i_{2} \in \mathcal{I}$ and $i_{1} \preceq i_{2}$ imply $\theta_{i_{1}} \leq \theta_{i_{2}}$.
\ed
We denote the set of  isotonic vectors in $\mathbb{R}^{s}$, with $s \leq \infty$, by $\bm{\mathcal{F}}^{is} = \bm{\mathcal{F}}^{is}(\mathcal{I})$. In the case of an infinite index set we consider the square summable vectors (with weights $\bm{w}$) from $\bm{l}_{2}^{\bm{w}}$, the Hilbert space of all square summable vectors with weights $\bm{w}$.

\bd
A  vector $\bm{\theta}^{*}\in \mathbb{R}^{s}$, with $s \leq \infty$, is the isotonic regression of an arbitrary vector $\bm{\theta} \in \mathbb{R}^{s}$ (or $\bm{\theta} \in \bm{l}_{2}^{\bm{w}}$, if $s=\infty$) over the pre-ordered index set $\mathcal{I}$ with weights $\bm{w} \in \mathbb{R}^{s}_{+}$ if 
\begin{eqnarray*}
\bm{\theta}^{*} = \underset{\bm{\xi} \in \bm{\mathcal{F}}^{is}}{\argmin} \sum_{i \in \mathcal{I}}(\xi_{i} - \theta_{i})^{2}w_{i}.
\end{eqnarray*}
\ed

Given a set $\mathcal{X}$ with a pre-order $\preceq$ on it one can generate a pre-order on the set $\mathcal{I}=\{1, 2 \dots \}$ of indices of the domain in $\mathcal{X}$ as follows. For $i_{1},  i_{2} \in \mathcal{I}$, $i_{1} \preceq  i_{2}$ if and only if $x_{i_{1}} \preceq x_{i_{2}}$. This pre-order on the index set $\mathcal{I}$ will be called the pre-order induced by the set $\mathcal{X}$ and will be denoted by the same symbol $\preceq$. Conversely, if one starts with the set $\mathcal{I}$ consisting of the indices of the elements in $\mathcal{X}$, and $\preceq$ is a pre-order on $\mathcal{I}$, the above correspondence defines a pre-order on $\mathcal{X}$. Therefore, in the sequel of the paper a bold symbol, e.g. $\bm{g}$, will denote a vector in $\mathbb{R}^{s}$, with $s \leq \infty$, whose $i$-th component is given by $g_{i} = g(x_{i})$, for $i = 1, \dots, s$, where $g(x)$ is a bounded real valued function on $\mathcal{X}$. In this case we will say that the vector $\bm{g}$ corresponds to the function $g(x)$ on $\mathcal{X}$ and vice versa.

\bnote
A real valued function $f(x)$ on the countable set $\mathcal{X}$ with the pre-order $\preceq$, defined on it, is isotonic if and only if its corresponding vector $\bm{f} \in \mathbb{R}^{s}$, with $s \leq \infty$, is an isotonic vector with respect to the corresponding pre-order $\preceq$ on its index set $\mathcal{I} = \{1, 2, \dots\}$, induced by the pre-order on $\mathcal{X}$. A real valued function $g^{*}(x)$ on the set $\mathcal{X}$ is the isotonic regression of a function $g(x)$ with weights $\bm{w}$ if and only if its corresponding vector $\bm{g}^{*} \in \mathbb{R}^{s}$ is the isotonic regression of the vector $\bm{g} \in \mathbb{R}^{s}$ with respect to the corresponding pre-order $\preceq$ on its index set $\mathcal{I} = \{1, 2, \dots \}$ with weights $\bm{w}$.
\enote

To state the inference problem treated in this paper, suppose that $\mathcal{X}$ is a finite or an infinite countable pre-ordered set and $\mathring{g}\in\mathcal{F}^{is}$ is a fixed unknown function. Suppose we are given observations $z_{i}$, $i= 1, \dots, n$, independent or not, that depend on the (parameter) $\mathring{g}$ in some way. In the sequel we will treat in detail  two important cases:\\ 
The data $z_{1}, \dots, z_{n}$ are observations of either of
\begin{enumerate}[label=(\roman*)]
\item $Z_{i}$, $i = 1, \dots, n$ independent identically distributed random variables taking values in $\mathcal{X}$, with probability mass function $\mathring{g}$.
\item $Z_{i} = (\bm{x}_{i}, Y_{i})$, $i = 1, \dots, n$, with $\bm{x}_{i}$ deterministic (design) points in $\mathcal{X}$ and $Y_{i}$ real valued random variables defined in the regression model
\begin{eqnarray*}
Y_{i} = \mathring{g}(x_{i}) + \varepsilon_{i}, \quad i=1, \dots, n,
\end{eqnarray*}
where $\varepsilon_{i}$ is a sequence of identically distributed random variables with $\mathbb{E}[\varepsilon_{i}] = 0$ and for simplicity we assume that $\mathrm{Var}[\varepsilon_{i}] = \sigma^{2} < \infty$. Note, that we will assume constant variance merely for simplicity. In fact, a heteroscedastic model would give us qualitatively similar results.
\end{enumerate}

Now assume that 
\begin{eqnarray*}
\hat{\bm{g}}_{n} = \hat{\bm{g}}_{n}(z_{1}, \dots, z_{n})
\end{eqnarray*}
is a $\mathbb{R}^{s}$-valued statistic. We will call the sequence $\{\hat{\bm{g}}_{n}\}_{n \geq 1}$ the basic estimator of $\mathring{\bm{g}}$. In order to discuss consistency and asymptotic distribution results we introduce the following basic topologies: When $s < \infty$, we study the Hilbert space with the inner product  $\langle \bm{g}_{1},\bm{g}_{2} \rangle = \sum_{i=1}^{s}g_{1, i}g_{2, i}w_{i}$, for $\bm{g}_{1}, \bm{g}_{2} \in \mathbb{R}^{s}$, endowed with its Borel $\sigma$-algebra $\mathcal{B} = \mathcal{B}(\mathbb{R}^{s})$ and when $s=\infty$ we study the space $\bm{l}_{2}^{\bm{w}}$ with the inner product $\langle \bm{g}_{1},\bm{g}_{2} \rangle = \sum_{i=1}^{\infty}g_{1, i}g_{2, i}w_{i}$, for a fixed weight vector $\bm{w}$ satisfying 
\begin{eqnarray}\label{conditw}
\begin{cases}
\underset{i}{\inf} \{ w_{i} \} >  0\\
\underset{i}{\sup} \{ w_{i} \} <  \infty,
\end{cases}
\end{eqnarray}
and we equip $\bm{l}_{2}^{\bm{w}}$ with its Borel $\sigma$-algebra $\mathcal{B} = \mathcal{B}(\bm{l}_{2}^{\bm{w}})$. 

Now define the isotonized estimator $\hat{\bm{g}}^{*}_{n}$ by 
\begin{eqnarray}\label{isotestg}
\hat{\bm{g}}^{*}_{n} = \underset{\bm{\xi} \in \bm{\mathcal{F}}^{is}}{\argmin} \sum_{i \in \mathcal{I}}(\xi_{i} - \hat{g}_{n,i})^{2}w_{i}.
\end{eqnarray}
In Section 5 we provide an example when the weights are not all the same and depend on data. The main goal of this paper is to study the asymptotic behaviour of $\hat{\bm{g}}^{*}_{n}$, as $n\to\infty$.

We make the following assumptions on the basic estimator $\hat{\bm{g}}_{n}$, for the finite, $s < \infty$, and the infinite, $s = \infty$, support case, respectively.\\
\bassum\label{astheta}
Suppose that $s < \infty$. Assume that $\hat{\bm{g}}_{n} \stackrel{p}{\to} \mathring{\bm{g}}$ for some $\mathring{\bm{g}} \in \bm{\mathcal{F}}^{is}$ and $B_{n}(\hat{\bm{g}}_{n} - \mathring{\bm{g}}) \stackrel{d}{\to} \bm{\lambda}$, where $\bm{\lambda}$ is a random vector in $(\mathbb{R}^{s}, \mathcal{B})$ and  $B_{n}$ is a diagonal $s\times s$ matrix with elements $[B_{n}]_{i i} = n^{q_{i}}$ with $q_{i}$ being real positive numbers.
\eassum
\bassum\label{asthetainf}
Suppose that $s = \infty$. Assume that $\hat{\bm{g}}_{n} \stackrel{p}{\to} \mathring{\bm{g}}$ for some $\mathring{\bm{g}} \in \bm{\mathcal{F}}^{is}$, and $B_{n}(\hat{\bm{g}}_{n} - \mathring{\bm{g}}) \stackrel{d}{\to} \bm{\lambda}$, where $\bm{\lambda}$ is a random vector in $(\bm{l}_{2}^{\bm{w}},\mathcal{B})$ and $B_{n}$ is a linear operator $\bm{l}_{2}^{\bm{w}} \to \bm{l}_{2}^{\bm{w}}$, such that for any $\bm{g} \in \bm{l}_{2}^{\bm{w}}$ it holds that $(B_{n}\bm{g})_{i} = n^{q_{i}} g_{i}$, with $q_{i}$ being the real positive numbers.
\eassum
Note that the matrix $B_{n}$ in Assumption \ref{astheta} and the operator $B_{n}$ in Assumption \ref{asthetainf}  allow for different rates of convergence for different components of $\hat{\bm{g}}_{n}$, i.e. the rates $q_{i}$ can be all the same but they do need to. The values of $q_{i}$ will be specified later. 

\section{The case of finitely supported functions}\label{sec:finitecase}
Let us assume that $s < \infty$, i.e. that the basic estimator  $\{\hat{\bm{g}}_{n}\}_{n \geq 1}$ is a sequence of finite-dimensional vectors. The next lemma states some well-known general properties of the isotonic regression of a finitely supported function.
\blem\label{propisot}
Suppose Assumption \ref{astheta} holds. Let $\hat{\bm{g}}^{*}_{n} \in \mathbb{R}^{s}$ be the isotonic regression of the vector $\hat{\bm{g}}_{n}$, defined in (\ref{isotestg}), for $n = 1, 2, 3, \dots$. Then the following hold:
\begin{enumerate}[label=(\roman*)]
\item $\hat{\bm{g}}^{*}_{n}$ exists and it is unique.
\item $\sum_{i=1}^{s} \hat{g}_{n,i} w_{i} = \sum_{i=1}^{s} \hat{g}^{*}_{n,i} w_{i}$, for all $n = 1, 2, \dots$. 
\item $\hat{\bm{g}}^{*}_{n}$, viewed as a mapping from $\mathbb{R}^{s}$ into $\mathbb{R}^{s}$, is continuous. Moreover, it is also continuous if it is viewed as a function on the $2s$-tuples of real numbers $(w_{1}, w_{2}, \dots, w_{s}, g_{1}, g_{2}, \dots, g_{s})$, with $w_{i} > 0$. 
\item  Assume that $a \leq \hat{g}_{n,i} \leq b$ holds for some constants $-\infty <a < b < \infty$, for all $n = 1, 2, \dots $ and $i = 1, \dots, s$. Then $\hat{\bm{g}}^{*}_{n}$ satisfies the same bounds as the basic estimator, i.e. $a \leq \hat{g}^{*}_{n, i} \leq b$, for all $n = 1, 2, \dots $ and $i = 1, \dots, s$. 
\item If $\hat{\bm{g}}_{n} \stackrel{p}{\to} \mathring{\bm{g}}$, then $\hat{\bm{g}}^{*}_{n} \stackrel{p}{\to} \mathring{\bm{g}}$.
\item $(\hat{\bm{g}}_{n} + \bm{c})^{*} = \hat{\bm{g}}^{*}_{n} + \bm{c}$ for all constant vectors $\bm{c} \in \mathbb{R}^{s}$,\\
$(c\hat{\bm{g}}_{n})^{*} = c\hat{\bm{g}}^{*}_{n} $ for all $c \in \mathbb{R}_{+}$.
\end{enumerate}  
\elem

We make a partition of the original set into comparable sets $\mathcal{X}^{(1)}, \dots, \mathcal{X}^{(k)}$
\begin{eqnarray}\label{partX}
\mathcal{X} = \cup_{v=1}^{k} \mathcal{X}^{(v)},
\end{eqnarray}
where each partition set $\mathcal{X}^{(v)}$ contains elements such that if $x \in \mathcal{X}^{(v)}$, then $x$ is comparable with at least one different element in $\mathcal{X}^{(v)}$ (if there are any), but not with any other element in $\mathcal{X}^{(\mu)}$ for any $\mu \neq v$.  In fact, the partition can be constructed even for an infinite set $\mathcal{X}=\{x_1,x_2,\ldots\}$, since a pre-order can be represented as the directed graph and any graph can be partitioned into isolated connected components and the partition is unique, cf. \cite{DuPardalos}.

Now assume that $\mathcal{X}$ is finite, that we are given the partition (\ref{partX}) and let $g^{(v)}(x)$, for $v=1, \dots, k$, be real valued functions defined on the sets $\mathcal{X}^{(v)}$ as $g^{(v)}(x) = g(x)$, whenever $x \in \mathcal{X}^{(v)}$, i.e. $g^{(v)}(x)$ is the restriction of the function $g(x)$ to the set $\mathcal{X}^{(v)}$. The family of functions $g^{(v)}(x)$ defined on the set $\mathcal{X}^{(v)}$, which are isotonic with respect to the pre-order, will be denoted by $\mathcal{F}^{is}_{(v)}$, for $v = 1, \dots, k$.

The next lemma states a natural result of an isotonic regression on $\mathcal{X}$, that it can be obtained as a concatenation of the individual isotonic regressions of the restrictions to the comparable sets.
\blem\label{isotoutpartlem}
Let $g(x)$ be an arbitrary real valued function on the finite set $\mathcal{X}$ with a pre-order $\preceq$ defined on it, and assume that the partition (\ref{partX}) is given. Then the isotonic regression of $g(x)$ with any positive weights $\bm{w}$ with respect to the pre-order $\preceq$ is equal to
\begin{eqnarray}\label{isonpartX}
g^{*}(x) = g^{*(v)}(x), \text{ whenever $x \in \mathcal{X}^{(v)}$},
\end{eqnarray}
where $g^{*(v)}(x)$ is the isotonic regression of the function $g^{(v)}(x)$ over the set $\mathcal{X}^{(v)}$ with respect to the pre-order $\preceq$. 
\elem

Now let $\mathring{g}(x)$ be the fixed function defined in Assumption \ref{astheta}, assume that we are given the partition (\ref{partX}) of $\mathcal{X}$ and for an arbitrary but fixed $v \in \{1, \dots, k \}$ let $\mathring{g}_{v}(x)$ be the restriction of $\mathring{g}(x)$ to $\mathcal{X}^{(v)}$. Then if $N_{v} = |\mathcal{X}^{(v)}|$ we can introduce the vector $\mathring{\bm{g}}_{v} = (\mathring{g}_{v,1}, \dots, \mathring{g}_{v,N_{v}}) = (\mathring{g}_{v}(x_{i_{1}}), \dots, \mathring{g}_{v}(x_{i_{N_{v}}}))$, where $x_{i_{1}}, \dots, x_{i_{N_{v}}}$ are the unique points in $\mathcal{X}^{(v)}$. Given $\mathring{g}_{v}(x)$ we can partition the set $\mathcal{X}^{(v)}$ into $m_{v}$ sets
\begin{eqnarray}\label{partXv}
\mathcal{X}^{(v)} = \cup_{l=1}^{m_{v}}\mathcal{X}^{(v, l)}.
\end{eqnarray}

The partition is constructed in the following way: We note first that the $N_{v}$ values in the vector $\mathring{\bm{g}}_{v}$ are not necessarily all unique, so there are $\tilde{m}_{v} \leq N_{v}$ unique values in $\mathring{\bm{g}}_{v}$. Then in a first step we construct $\tilde{m}_{v}$ level sets
\begin{eqnarray*}\label{}
\tilde{\mathcal{X}}^{(v, l)} = \{x \in \mathcal{X}^{(v)}: \mathring{g}_{v}(x_{i}) =  \mathring{g}_{v,l}\}
\end{eqnarray*}
with $l = 1, \dots, \tilde{m}_{v}$.

Next we note that for any non-singleton level set $\tilde{\mathcal{X}}^{(v, l)}$ there might be non-comparable points, i.e. $x_{i}, x_{j} \in \tilde{\mathcal{X}}^{(v, l)}$ can be such that neither $x_{i} \preceq x_{j}$ nor $x_{j} \preceq x_{i}$ hold. Therefore, in the second step for each fixed $l$ we can partition (if necessary) the level set  $\tilde{\mathcal{X}}^{(v, l)}$ into sets with comparable elements, analogously to the construction of (\ref{partX}). We can do this for every $v$ and end up in a partition (\ref{partXv}) with $\tilde{m}_{v} \leq m_{v} \leq N_{v}$. 

In the partition (\ref{partXv}) each set $\mathcal{X}^{(v, l)}$ is characterised by
\begin{enumerate}[label=(\roman*)]
\item for every $x \in \mathcal{X}^{(v, l)}$ we have $\mathring{g}_{v}(x) = \mathring{g}_{v, l}$,
\item if $|\mathcal{X}^{(v, l)}| \geq 2$ then for every $x \in \mathcal{X}^{(v, l)}$ there is at least one $x' \in \mathcal{X}^{(v, l)}$ such that $x \sim x'$.
\end{enumerate}
We have therefore proved the following lemma.
\blem\label{corpartX}
For any countable set $\mathcal{X}$ with the pre-order $\preceq$ and any isotonic function $\mathring{g}(x)$, defined on it, there exists a unique partition $\mathcal{X} =  \cup_{v=1}^{k} \cup_{l=1}^{m_{v}}\mathcal{X}^{(v, l)}$, satisfying the statements $(i)$ and $(ii)$ above. For the index set $\mathcal{I}$ with the pre-order $\preceq$ generated by the set $\mathcal{X}$ and any isotonic function $\mathring{g}(x)$, defined on $\mathcal{X}$, there exists a unique partition $\mathcal{I} =  \cup_{v=1}^{k} \cup_{l=1}^{m_{v}}\mathcal{I}^{(v, l)}$, satisfying conditions analogous to $(i)$ and $(ii)$ stated above.
\elem
\bd\label{complsdef}
The set $\mathcal{X}$ will be called decomposable  if in the partition, defined in (\ref{partX}), $k > 1$. In the partition (\ref{partXv}) the sets $\mathcal{X}^{(v, l)}$ will be called the comparable level sets of $\mathring{g}(x)$. In the corresponding partition of the index set $\mathcal{I}$ the sets $\mathcal{I}^{(v, l)}$ will be called the comparable level index sets of $\mathring{\bm{g}}$.
\ed

Recall that $g^{(v,l)}(x)$ is the restriction of the function $g(x)$ to the comparable level set $\mathcal{X}^{(v, l)}$, for $l = 1, \dots, m_{v}$ and $v = 1, \dots, k$. 

In the case of a non-decomposable  set, the full partition will be written as $\mathcal{X} =  \cup_{l=1}^{m_{1}}\mathcal{X}^{(1, l)} \equiv \cup_{l=1}^{m}\mathcal{X}^{(l)}$, so we may then drop the index $v=1$. Similarly, in this case $g^{(1,l)}(x) \equiv g^{(l)}(x)$ denotes the restriction of a function $g(x)$ to the comparable level set $\mathcal{X}^{(l)} \equiv \mathcal{X}^{(1, l)}$.

Next, suppose that $\mathcal{X}$ is a non-decomposable set, and let us consider an arbitrary function $\mathring{g}(x) \in \mathcal{F}^{is}$.  Assume that for $\mathring{g}(x)$ there has been made a partition $\mathcal{X} = \cup_{l=1}^{m}\mathcal{X}^{(l)}$ in (\ref{partXv}), satisfying $(i)$ and $(ii)$. Define the smallest comparable level distance of $\mathring{g}$ as 
\begin{eqnarray}\label{varepstild}
\tilde{\varepsilon} &=& \inf \{|\mathring{g}_{l'} - \mathring{g}_{l}|: l, l' = 1, \dots, m, \,  l \neq l', \\\nonumber
&& \exists x_{1} \in \mathcal{X}^{(l)}, \exists x_{2} \in \mathcal{X}^{(l')}, \text{ such that } x_{1} \sim x_{2} \},
\end{eqnarray}
Note, that $\tilde{\varepsilon}$ is always finite and for the finite support case, $s < \infty$, also $\tilde{\varepsilon} > 0$.

\blem\label{isonlevels}
Consider an arbitrary real valued function $g(x)$ on a non-\\decomposable finite set $\mathcal{X}$ with the pre-order $\preceq$ and let $\tilde{\varepsilon}$ be defined in (\ref{varepstild}). If 
\begin{eqnarray}\label{levconditfin}
\underset{x \in \mathcal{X}}{\sup} \{ |g(x) - \mathring{g}(x)|\} < \tilde{\varepsilon}/2,
\end{eqnarray}
then the isotonic regression of $g(x)$ is given by
\begin{eqnarray}\label{isregndcs}
g^{*}(x) = g^{*(l)}(x), \text{ whenever $x \in \mathcal{X}^{(l)}$},
\end{eqnarray}
where $g^{*(l)}(x)$ is the isotonic regression of the function $g^{(l)}(x)$ over the set $\mathcal{X}^{(l)}$ with respect to the pre-order $\preceq$. Therefore, the function $g^{*}(x)$ is a concatenation of the isotonic regressions of the restrictions of $g(x)$ to the comparable level sets of $\mathring{g}(x)$.
\elem

The next lemma is an auxiliary result which will be used later in the proof of the asymptotic distribution of $\hat{\bm{g}}_{n}^{*}$.
\blem\label{propconv}
Assume $\bm{X}_{n}$ and $\bm{Y}_{n}$ are sequences of random vectors, taking values in the space $\mathbb{R}^{s}$, for $s \leq \infty$, with some metric on it, endowed with its Borel $\sigma$-algebra. If $\bm{X}_{n} \stackrel{d}{\to} \bm{X}$ and $\lim_{n\to\infty}\mathbb{P}[\bm{X}_{n} = \bm{Y}_{n}] = 1$, then $\bm{Y}_{n} \stackrel{d}{\to} \bm{X}$. 
\elem

Let us consider the sequence $B_{n}(\hat{\bm{g}}_{n}^{*} - \mathring{\bm{g}})$, where $\hat{\bm{g}}_{n}^{*}$ is the isotonic regression of $\hat{\bm{g}}_{n}$, which was defined in Assumption \ref{astheta}, and with a specified matrix $B_{n}$. As mentioned in Assumption \ref{astheta}, we allow different rates of convergence $n^{q_{i}}$ for different components of $\hat{\bm{g}}_{n}$. We however require $q_{i}$, for $i = 1, \dots, s$, to be equal on the comparable level index sets $\mathcal{I}^{(v, l)}$ of $\mathring{\bm{g}}$, i.e. $q_{i}$, for $i = 1, \dots, s$, are real positive numbers such that $q_{i_{1}} = q_{i_{2}}$, whenever $i_{1}, i_{2} \in \mathcal{I}^{(v, l)}$. 

We introduce an operator $\varphi: \mathbb{R}^{s} \to \mathbb{R}^{s}$ defined in the following way. First, for any vector $\bm{\theta} \in \mathbb{R}^{s}$ we define the coordinate evaluation map $\theta(x): \mathcal{X} \to \mathbb{R}$, corresponding to the vector $\bm{\theta}$, by $\theta(x_{i}) = \theta_{i}$, for $i = 1, \dots, s$. Then, let $\theta^{*(v',l')}(x)$ be the isotonic regression of the restriction of $\theta(x)$ to the comparable level set $\mathcal{X}^{(v', l')}$ of $\mathring{g}(x)$, and define
\begin{eqnarray}\label{phiop}
\varphi(\bm{\theta})_{i} = \theta^{*(v',l')}(x_{i}),
\end{eqnarray} 
for $i=1, \dots, s$, with $(v',l')$ the (unique) indices such that $x_{i} \in \mathcal{X}^{(v', l')}$.

The asymptotic distribution of $B_{n}(\hat{\bm{g}}_{n}^{*} - \mathring{\bm{g}})$ is given in the following theorem.  
\bth\label{thmasymfin}
Suppose that Assumption \ref{astheta} holds. Then 
\begin{eqnarray}\label{asisot}
B_{n}(\hat{\bm{g}}_{n}^{*} - \mathring{\bm{g}}) \stackrel{d}{\to} \varphi(\bm{\lambda}),
\end{eqnarray}
where $\varphi$ is the operator, defined in (\ref{phiop}). \eth
\prf

First, from Lemma \ref{corpartX} we have that any pre-ordered set $\mathcal{X}$ can be uniquely partitioned as 
\begin{eqnarray}\label{partXsumm}
\mathcal{X} &=& \cup_{v=1}^{k} \mathcal{X}^{(v)},\nonumber\\
\mathcal{X}^{(v)} &=& \cup_{l=1}^{m_{v}}\mathcal{X}^{(v, l)},
\end{eqnarray}
and with the partition (\ref{partXsumm}) of $\mathcal{X}^{(v)}$ determined by the isotonic vector $\mathring{\bm{g}}$. 

Second, as shown in Lemma \ref{isotoutpartlem}, the isotonic regression of $g(x)$ on the original set $\mathcal{X}$ can be obtained as a concatenation of the separate isotonic regressions of the restrictions of $g(x)$ to the non-decomposable sets in the partition (\ref{partX}). Therefore, without loss of generality, we can assume that the original set $\mathcal{X}$ is non-decomposable. Thus, any $x \in \mathcal{X}$ is comparable with at least one different element of $\mathcal{X}$, $k=1$, and 
\begin{eqnarray*}
\mathcal{X} &=& \cup_{l=1}^{m_{1}} \mathcal{X}^{(1, l)}\\
&\equiv& \cup_{l=1}^{m} \mathcal{X}^{(l)}
\end{eqnarray*}
and $\mathring{g}_{1, l} \equiv \mathring{g}_{l}$. Note, that we have dropped the index $v$. 

Third, since $\hat{\bm{g}}_{n}$ is consistent, by Assumption \ref{astheta}, for any $\varepsilon>0$, 
\begin{eqnarray}\label{probdiffg}
\mathbb{P}[\underset{x \in \mathcal{X}}{\sup}\{ |\hat{g}_{n}(x) - \mathring{g}(x)| \} < \varepsilon]\to 1,
\end{eqnarray}
as $n\to\infty$. Note that the comparable level distance $\tilde{\varepsilon}$ of $\mathring{g}$, defined in  (\ref{varepstild}), satisfies $\tilde{\varepsilon} > 0$, and take $\varepsilon = \tilde{\varepsilon}/2$. Then from Lemma \ref{isonlevels} we obtain
\begin{eqnarray}\label{equivepsphi}
\{ \underset{x \in \mathcal{X}}{\sup}\{ |\hat{g}_{n}(x) - \mathring{g}(x)| \} < \epsilon\} \subseteq 
\{\hat{\bm{g}}_{n}^{*} = \varphi(\hat{\bm{g}}_{n}) \}.
\end{eqnarray}
Therefore, (\ref{probdiffg}) and (\ref{equivepsphi}) imply 
\begin{eqnarray}\label{Pequivepsphi}
\mathbb{P}[\hat{\bm{g}}_{n}^{*} = \varphi(\hat{\bm{g}}_{n})] \to 1,
\end{eqnarray}
as $n \to \infty$.

Next, since the isotonic regression is a continuous map (statement $(iii)$ of Lemma \ref{propisot}), the operator $\varphi$ is a continuous map from $\mathbb{R}^{s}$ to $\mathbb{R}^{s}$. Therefore, using the continuous mapping theorem, cf. \cite{vandervaart:1998}, we get
\begin{eqnarray}\label{convphiBng}
\varphi(B_{n}(\hat{\bm{g}}_{n} - \mathring{\bm{g}})) \stackrel{d}{\to} \varphi(\bm{\lambda}).
\end{eqnarray}

Furthermore, using statement (vi) of Lemma \ref{propisot} and taking into account the definition of the matrix $B_{n}$, we get
\begin{eqnarray}\label{equalpBBp}
\varphi(B_{n}(\hat{\bm{g}}_{n} - \mathring{\bm{g}})) = B_{n}(\varphi(\hat{\bm{g}}_{n}) - \mathring{\bm{g}}).
\end{eqnarray}

Then (\ref{Pequivepsphi}), (\ref{convphiBng}) and (\ref{equalpBBp}) imply that 
\begin{eqnarray}\label{Pequivexprisot}
\mathbb{P}[B_{n}(\hat{\bm{g}}_{n}^{*} - \mathring{\bm{g}}) = B_{n}(\varphi(\hat{\bm{g}}_{n}) - \mathring{\bm{g}})] \to 1,
\end{eqnarray}  
as $n\to\infty$. Finally, using Lemma \ref{propconv}, from (\ref{convphiBng}) and (\ref{Pequivexprisot}) we prove that
\begin{eqnarray*}
B_{n}(\hat{\bm{g}}_{n}^{*} - \mathring{\bm{g}}) \stackrel{d}{\to} \varphi(\bm{\lambda}),
\end{eqnarray*}
as $n\to\infty$.
\eop

For a given pre-order $\preceq$ on $\mathcal{X}$, the statement that there exists a matrix $\bm{A}$ such that $\bm{A}\bm{g} \geq \bm{0}$ is equivalent to the statement that $\bm{g}$ is isotonic with respect to $\preceq$, cf. Proposition 2.3.1 in \cite{silvapsen}. Therefore, if there are no linear constraints imposed on the basic estimator $\hat{\bm{g}}_{n}$, Theorem \ref{thmasymfin} can also be established by using the results on estimation when a parameter is on a boundary, in Section 6 in \cite{Andrews}. 

Assume that each vector $\hat{\bm{g}}_{n}$ has the following linear constraint $\sum_{i=1}^{s} \hat{g}_{n,i} w_{i} \\= c$ (for example, in the case of estimation of a probability mass function it would be $\sum_{i=1}^{s} \hat{g}_{n,i} = 1$). 
Then, the expression for a limiting distribution in Theorem \ref{thmasymfin} does not follow directly from the results in \cite{Andrews} in the case when $\hat{\bm{g}}_{n}$ is linearly constrained. However, the result of Theorem \ref{thmasymfin} holds, because, as established in statement $(ii)$ of Lemma \ref{propisot}, isotonic regression with weights $\bm{w}$ preserves the corresponding linear constraint.

Also, we note that a result similar to Theorem \ref{thmasymfin} for the isotonic regression of a finitely supported empirical pmf, was obtained in the proof of Theorem 5.2.1 in \cite{robertsonorder}.

Next we consider the case when the vector of weights $\bm{w}$ is not a constant, i.e. we assume that some non-random sequence $\{\bm{w}_{n}\}_{n\geq1}$, where each vector $\bm{w}_{n}$ satisfies the condition (\ref{conditw}), converges to some non-random vector $\bm{w}$, which also satisfies (\ref{conditw}). We denote by $\theta^{*\bm{w}}(x)$ the isotonic regression of $\theta(x)$ with weights $\bm{w}$ and  analogously to (\ref{phiop}) we introduce the notation $\varphi^{\bm{w}}(\bm{\theta})$
\begin{eqnarray}\label{phiopw}
\varphi^{\bm{w}}(\bm{\theta})_{i} = \theta^{*\bm{w}(v',l')}(x_{i}),
\end{eqnarray} 
where $\theta^{*\bm{w}(v',l')}(x)$ is the isotonic regression, with weights $\bm{w}$, of the restriction of $\theta(x)$ to the comparable level set $\mathcal{X}^{(v', l')}$ of $ \mathring{g}(x)$, where the indices $v'$ and $l'$ are such that $x_{i} \in \mathcal{X}^{(v', l')}$. Define the isotonic regression $\hat{\bm{g}}_{n}^{*\bm{w}_{n}}$ of the basic estimator $\hat{\bm{g}}_{n}$. The next theorem gives the limiting distribution of $\hat{\bm{g}}_{n}^{*\bm{w_{n}}}$.

\bth\label{thmasymfinw}
Suppose that Assumption \ref{astheta} holds. Then the asymptotic distribution of the isotonic regression $\hat{\bm{g}}_{n}^{*\bm{w_{n}}}$ of the basic estimator $\hat{\bm{g}}_{n}$ is given by
\begin{eqnarray}\label{asisot}
B_{n}(\hat{\bm{g}}_{n}^{*\bm{w}_{n}} - \mathring{\bm{g}}) \stackrel{d}{\to} \varphi^{\bm{w}}(\bm{\lambda}),
\end{eqnarray}
where $\varphi^{\bm{w}}$ is the operator, defined in (\ref{phiopw}). \eth
\prf
Without loss of generality, we can assume that the original set $\mathcal{X}$ is non-decomposable. First, since the sequence $\hat{\bm{g}}_{n}$ is consistent, then for any  
\begin{eqnarray}\label{probdiffgw}
\mathbb{P}[\underset{x \in \mathcal{X}}{\sup}\{ |\hat{g}_{n}(x) - \mathring{g}(x)| \} < \tilde{\varepsilon}/2]\to 1,
\end{eqnarray}
as $n\to\infty$, with $\tilde{\varepsilon}$ defined in \ref{varepstild}. Using the statement of Lemma \ref{isonlevels}, we obtain
\begin{eqnarray}\label{equivepsphiw}
\{ \underset{x \in \mathcal{X}}{\sup}\{ |\hat{g}_{n}(x) - \mathring{g}(x)| \} < \tilde{\varepsilon}/2\} \subseteq 
\{\hat{\bm{g}}_{n}^{*\bm{w}_{n}} = \varphi^{\bm{w}_{n}}(\hat{\bm{g}}_{n}) \}.
\end{eqnarray}
Note that the result of Lemma \ref{isonlevels} holds for any weights $\bm{w}_{n}$.

Therefore, from (\ref{probdiffgw}) and (\ref{equivepsphiw}) we have 
\begin{eqnarray}\label{Pequivepsphiw}
\mathbb{P}[\hat{\bm{g}}_{n}^{*\bm{w}_{n}} = \varphi(\hat{\bm{g}}_{n})] \to 1,
\end{eqnarray}
as $n \to \infty$.

Second, from statement $(iii)$ of Lemma \ref{propisot}, the operators $\varphi^{\bm{w}_{n}}$, $\varphi^{\bm{w}}$ are continuous maps from $\mathbb{R}^{2s}$ to $\mathbb{R}^{s}$, for all weights $\bm{w}_{n}$, $\bm{w}$ satisfying (\ref{conditw}). Using the (extended) continuous mapping theorem, cf. \cite{vandervaart:1998}, we get
\begin{eqnarray}\label{convphiBngw}
\varphi^{\bm{w}_{n}}(B_{n}(\hat{\bm{g}}_{n} - \mathring{\bm{g}})) \stackrel{d}{\to} \varphi^{\bm{w}}(\bm{\lambda}),
\end{eqnarray}
where $\bm{w}$ is the limit of the sequence $\{\bm{w}_{n}\}_{n\geq1}$.

Third, using statement (vi) of Lemma \ref{propisot} and the definition of the matrix $B_{n}$ we obtain
\begin{eqnarray}\label{equalpBBpw}
\varphi^{\bm{w}_{n}}(B_{n}(\hat{\bm{g}}_{n} - \mathring{\bm{g}})) = B_{n}(\varphi^{\bm{w}_{n}}(\hat{\bm{g}}_{n}) - \mathring{\bm{g}}).
\end{eqnarray}
Therefore, (\ref{equalpBBpw}) gives us 
\begin{eqnarray}\label{Pequivexprisotw}
\mathbb{P}[B_{n}(\hat{\bm{g}}_{n}^{*\bm{w}_{n}} - \mathring{\bm{g}}) = B_{n}(\varphi^{\bm{w}_{n}}(\hat{\bm{g}}_{n}) - \mathring{\bm{g}})] \to 1,
\end{eqnarray}  
as $n\to\infty$. Finally, using Lemma \ref{propconv}, from (\ref{convphiBngw}) and (\ref{Pequivexprisotw}) we prove that
\begin{eqnarray*}
B_{n}(\hat{\bm{g}}_{n}^{*\bm{w}_{n}} - \mathring{\bm{g}}) \stackrel{d}{\to} \varphi^{\bm{w}}(\bm{\lambda}),
\end{eqnarray*}
as $n\to\infty$.
\eop

\section{The case of infinitely supported functions}\label{sec:infinitecase}
In this section we assume that the original set $\mathcal{X} = \{ x_{1}, x_{2} \dots \}$ is an infinite countable enumerated  set with a pre-order $\preceq$ defined on it. 

In the case of infinitely supported functions the isotonic regression's properties are similar to the ones in the finite case, but the proofs are slightly different. For completeness we state these properties in the following lemma.
\blem\label{propisotinf}
Suppose Assumption \ref{astheta} holds. Let $\hat{\bm{g}}^{*}_{n} \in \bm{l}_{2}^{\bm{w}}$ be the isotonic regression of the vector $\hat{\bm{g}}_{n} \in \bm{l}_{2}^{\bm{w}}$, for $n = 1, 2, 3, \dots$. Then statements $(i)$ - $(vi)$ of Lemma \ref{propisot} hold, with $(iii)$ suitably changed to the mapping from $\bm{l}_{2}^{\bm{w}}$ to $\bm{l}_{2}^{\bm{w}}$.
\elem

We partition the original set $\mathcal{X}$ in the same way as it was done in the finite case, i.e., first, let
\begin{eqnarray}\label{partXinf}
\mathcal{X} = \cup_{v=1}^{k} \mathcal{X}^{(v)},
\end{eqnarray}
where $k \leq \infty$ is the number of sets and each set $\mathcal{X}^{(v)}$ is such that if $x \in \mathcal{X}^{(v)}$, then $x$ is comparable with at least one different element in $\mathcal{X}^{(v)}$ (if there are any), but not with any other elements which belong to other sets in the partition. Note, that it is possible for $\mathcal{X}^{(v)}$ to have only one element. The partition of $\mathcal{X}$ is unique and $\mathcal{X}^{(v)}\cap\mathcal{X}^{(v')}=\emptyset$, for $v\neq v'$. 

Furthermore, for the fixed function $\mathring{g} \in \mathcal{F}^{is}$, defined in Assumption \ref{asthetainf} (or, equivalently, its corresponding isotonic vector $\mathring{\bm{g}} \in \bm{\mathcal{F}}^{is}$), we partition each set $\mathcal{X}^{(v)}$ in (\ref{partXinf}) into the comparable level sets of $\mathring{g}$, i.e.
\begin{eqnarray}\label{partXvinf}
\mathcal{X}^{(v)} = \cup_{l=1}^{m_{v}}\mathcal{X}^{(v, l)},
\end{eqnarray}
in the same way as it was done in the finite case in (\ref{partXv}).

Note, that since $\mathring{\bm{g}} \in \bm{l}_{2}^{\bm{w}}$ and condition (\ref{conditw}) is satisfied, the cardinality of any set $\mathcal{X}^{(v, l)}$ is less than infinity whenever $\mathring{g}_{v,l} \neq 0$, otherwise we would have $\sum_{i=1}^{\infty}(\mathring{g}_{i})^{2}w_{i} = \infty$, which would mean that $\mathring{g} \not\in l_{2}^{\bm{w}}$. The set $\mathcal{X}^{(v, l)}$ can have infinitely many elements only if $\mathring{g}_{v,l} = 0$.

For the partition in (\ref{partXinf}) we obtain a result similar to the one obtained in Lemma \ref{isotoutpartlem} for the finite case.
\blem\label{propisotinf}
Let $g(x)$ be an arbitrary real valued function in $l_{2}^{\bm{w}}$ on the set $\mathcal{X}$ with a pre-order $\preceq$ defined on it. Then the isotonic regression of $g(x)$ with any positive weights $\bm{w}$ is equal to
\begin{eqnarray}\label{isonpartX}
g^{*}(x) = g^{*(v)}(x), \text{ whenever $x \in \mathcal{X}^{(v)}$},
\end{eqnarray}
where $g^{*(v)}(x)$ is the isotonic regression of the restriction of the function $g(x)$ to the set $\mathcal{X}^{(v)}$ over this set with respect to the pre-order $\preceq$. 
\elem

As a consequence of Lemma \ref{propisotinf}, without loss of generality in the sequel of the paper we can assume that the original set $\mathcal{X}$ is non-decomposable  and use the same notations as in the finite case, i.e. $\mathcal{X} = \cup_{l=1}^{m}\mathcal{X}^{(l)} \equiv \cup_{l=1}^{m_{1}}\mathcal{X}^{(1, l)}$ and, respectively, $g^{(l)}(x) \equiv g^{(1,l)}(x)$ for the restriction of the function $g(x)$ to the set $\mathcal{X}^{(l)}$.

In the case of an infinite support the result of Lemma \ref{isonlevels} is generally not applicable, because the value of $\tilde{\varepsilon}$ can in this case be zero. We therefore make the following slight modification of  Lemma \ref{isonlevels}. Thus, assume that for a function $\mathring{g}(x) \in \mathcal{F}^{is}$ we have made a partition $\mathcal{X} = \cup_{l=1}^{m}\mathcal{X}^{(l)}$ with $m \leq \infty$. Furthermore, for any finite positive integer number $m' < m \leq \infty$ we choose $m'$ comparable level sets $\mathcal{X}^{(l_{j})}$, such that the values of the function $\mathring{g}(x)$ on them satisfy $|\mathring{g}_{l_{1}}| \geq |\mathring{g}_{l_{2}}| \geq \dots \geq |\mathring{g}_{l_{m'}}|$. Next, we rewrite the partition as 
\begin{eqnarray}\label{partpartit}
\mathcal{X} = \mathcal{X}^{(l_{1})}\cup\mathcal{X}^{(l_{2})}\cup\dots\cup\mathcal{X}^{(l_{m'})}\cup\mathcal{X}^{(l_{m'+1})},
\end{eqnarray}
where $\mathcal{X}^{(l_{m'+1})} = \mathcal{X} \setminus  \mathcal{X}^{(l_{1})}\cup\mathcal{X}^{(l_{2})}\cup\dots\cup\mathcal{X}^{(l_{m'})}$. Define 
\begin{eqnarray}\label{vareptdpr}
\tilde{\varepsilon}' &=& \inf \{|\mathring{g}_{l'} - \mathring{g}_{l}|: l' \in \{l_{1}, \dots, l_{m'}\},\\\nonumber
&& l \in \{ 1, \dots, m \}, \exists x_{1} \in \mathcal{X}^{(l)},\\\nonumber
&& \exists x_{2} \in \mathcal{X}^{(l')}, \text{ such that } x_{1} \sim x_{2} \},
\end{eqnarray}
and note that $\tilde{\varepsilon}'$ is always positive.
\blem\label{isonlevelsinf}
Consider an arbitrary real valued function $g(x) \in l_{2}^{\bm{w}}$ on a non-decomposable  infinite countable set $\mathcal{X}$ with the pre-order $\preceq$ defined on it. Suppose that $\tilde{\varepsilon}'$ is defined in (\ref{vareptdpr}). If for some $\mathring{g}(x) \in \mathcal{F}^{is}$ we have
\begin{eqnarray*}
\underset{x \in \mathcal{X}}{\sup} \{ |g(x) - \mathring{g}(x)|\} < \tilde{\varepsilon}'/2, 
\end{eqnarray*}
then the isotonic regression of $g(x)$ is given by
\begin{eqnarray*}
g^{*}(x) = g^{*(l'')}(x), \text{ whenever $x \in \mathcal{X}^{(l'')}$ , for $l'' \in \{l_{1}, \dots, l_{m'}, l_{m' +1}\}$},
\end{eqnarray*}
where $g^{*(l'')}(x)$ is the isotonic regression of the function $g^{(l')}(x)$ over the set $\mathcal{X}^{(l')}$ with respect to the pre-order $\preceq$. Therefore, the function $g^{*}(x)$ is a concatenation of the isotonic regressions of the restrictions of $g(x)$ to the sets $\mathcal{X}^{(l_{1})},\mathcal{X}^{(l_{2})},\dots,\mathcal{X}^{(l_{m'})}$ and $\mathcal{X}^{(l_{m'+1})}$.
\elem

Next we state and prove an auxiliary lemma, see also Problem III.6.3 in \cite{Shrv2007}, which will be used in the final theorem.
\blem\label{convl2}
Let $\bm{Z}_{n}$, for $n = 1, \dots, \infty$, be a tight sequence of random vectors in $\bm{l}_{2}^{\bm{w}}$, endowed with its Borel $\sigma$-algebra $\mathcal{B}$. Consider the set of indices $\mathcal{I} = \{1, 2, \dots, \infty\}$ of the components of the vectors $\bm{Z}_{n}$. Assume that for some random vector $\bm{Z}$ in $(\bm{l}_{2}^{\bm{w}}, \mathcal{B})$ and some rearrangement $\tilde{\mathcal{I}}$ of the original index set $\mathcal{I}$ the following holds: For any positive finite integer $s$ we have $\bm{\tilde{Z}}^{(1, s)}_{n} \stackrel{d}{\to} \bm{\tilde{Z}}^{(1, s)}$, where $\bm{\tilde{Z}}^{(1, s)}_{n}$ and $\bm{\tilde{Z}}^{(1, s)}$ are vectors in $\mathbb{R}^{s}$  constructed from the elements of the vectors $\bm{Z}_{n}$ and $\bm{Z}$ in such a way that the $j$-th elements of $\bm{\tilde{Z}}^{(1, s)}_{n}$ and $\bm{\tilde{Z}}^{(1, s)}$ are equal to the $\tilde{i}_{j}$-th elements of the vectors $\bm{Z}_{n}$ and $\bm{Z}$, respectively, with $\tilde{i}_{j}$ being the $j$-th index from the rearranged index set $\tilde{\mathcal{I}}$. Then $\bm{Z}_{n} \stackrel{d}{\to} \bm{Z}$.
\elem

Finally, the next theorem gives the limiting distribution of $\hat{\bm{g}}_{n}^{*}$. Similarly to the finite case we introduce the operator $\varphi: \bm{l}_{2}^{\bm{w}} \to \bm{l}_{2}^{\bm{w}}$, defined in the following way. For any vector $\bm{\theta} \in \bm{l}_{2}^{\bm{w}}$ we consider the coordinate evaluation map $\theta(x): \mathcal{X} \to \mathbb{R}$ defined as $\theta(x_{i}) = \theta_{i}$, for $i = 1, \dots, \infty$. Then, let 
\begin{eqnarray}\label{phiinf}
\varphi(\bm{\theta})_{i} = \theta^{*(v',l')}(x_{i}),
\end{eqnarray}
where $\theta^{*(v',l')}(x)$ is the isotonic regression of the restriction of $\theta(x)$ to the set $\mathcal{X}^{(v', l')}$ in the partition of $\mathcal{X}$. The indices $v'$ and $l'$ are such that $x_{i} \in \mathcal{X}^{(v', l')}$. The restriction of $\varphi(\bm{\theta})$ to the comparable index level set $\mathcal{I}^{(v, l)}$ will be denoted by $[\varphi(\bm{\theta})]^{(v, l)}$
\bth\label{thmasyminf}
Suppose the Assumption \ref{asthetainf} holds. Then the asymptotic distribution of the isotonized estimator $\hat{\bm{g}}_{n}^{*}$ is given by
\begin{eqnarray}\label{asisot}
B_{n}(\hat{\bm{g}}_{n}^{*} - \mathring{\bm{g}}) \stackrel{d}{\to} \varphi(\bm{\lambda}),
\end{eqnarray}
where $\varphi$ is the operator defined in (\ref{phiinf}). \eth
\prf
Let us consider the partition of the original set $\mathcal{X} = \cup_{l=1}^{m}\mathcal{X}^{(l)}$ made for the function $\mathring{g}(x)$. As it was shown above, the  cardinality $|\mathcal{X}^{(l)}|$ of each comparable level set in the partition must be less than infinity, unless $\mathring{g}_{l} = 0$, in which case it can have infinite cardinality.  Note that if the number of terms in the partition is less than infinity, i.e. $m < \infty$, then some terms (or just one) in the partition are such that the function $\mathring{g}(x)$ is equal to zero on them, i.e. $\mathring{g}_{l} = 0$. Therefore, in this case we can use the same approach as in the case of the finite set $\mathcal{X}$ (Lemma \ref{isonlevels}), because in this case the smallest comparable level distance $\tilde{\varepsilon}$, defined in (\ref{varepstild}), is greater than zero.

Therefore, further in the proof we assume that $m=\infty$ and write the partition as $\mathcal{X} = \cup_{l=1}^{\infty}\mathcal{X}^{(l)}$. First, for any positive integer $m' < \infty$ let us take $m'$ terms from the partition of $\mathcal{X}$ which satisfy $|\mathring{g}_{l_{1}}| \geq |\mathring{g}_{l_{2}}| \geq \dots \geq |\mathring{g}_{l_{m'}}|$. 

Second, since the sequence $\hat{\bm{g}}_{n}$ is consistent, then for any 
$\varepsilon > 0$ 
$$
\lim_{n\to\infty}\mathbb{P}[\underset{i \in \mathcal{I}}{\sup} \{ |\hat{g}_{n,i} - \mathring{g}_{i}|\} < \varepsilon] = 1.
$$
Therefore, letting $\varepsilon = \tilde{\varepsilon}'/2$, with $\tilde{\varepsilon}'$ defined in (\ref{vareptdpr}), by Lemma \ref{isonlevelsinf} we obtain that, for the isotonic regression $\hat{\bm{g}}_{n}^{*}$ of $\hat{\bm{g}}_{n}$
\begin{eqnarray}\label{infconcat}
\lim_{n\to\infty}\mathbb{P}[\hat{g}^{*}_{i} = \hat{g}^{*(l'')}_{i}] =1, \text{ whenever $i \in \mathcal{I}^{(l'')}$} ,\\\nonumber \text{ for $l'' \in \{l_{1}, \dots, l_{m'+1}\}$,}
\end{eqnarray}
where $\mathcal{I}^{(l')}$, for $l' \in \{l_{1}, \dots, l_{m'}\}$, are the comparable level sets and $\mathcal{I}^{(m'+1)}$ is the index set of $\mathcal{X}^{(m' +1)} = \mathcal{X} \setminus  \mathcal{X}^{(l_{1})}\cup\mathcal{X}^{(l_{2})}\cup\dots\cup\mathcal{X}^{(l_{m'})}$.

Third, let us introduce a linear operator $A^{(m')}: \bm{l}_{2}^{\bm{w}} \to \mathbb{R}^{s}$, with \\ $s = \sum_{l\in\{l_{1}, \dots  l_{m'}\}} |\mathcal{X}^{(l)}|$, such that for any $\bm{g} \in \bm{l}_{2}^{\bm{w}}$ the first $|\mathcal{X}^{(l_{1})}|$ elements of the vector $A^{(m')} \bm{g}$ are equal to ones taken from $\bm{g}$ whose indices are in $\mathcal{I}^{(l_{1})}$, the second $|\mathcal{X}^{(l_{2})}|$ elements are the ones from $\bm{g}$ whose indices are from $\mathcal{I}^{(l_{2})}$ and so on. Therefore, using the result in (\ref{infconcat}), the definition of $B_{n}$ and statement (vi) of Lemma \ref{propisot}, we have that
\begin{eqnarray}\label{eqaulAphiinf}
\lim_{n\to\infty}\mathbb{P}[A^{(m')} \varphi(B_{n}(\hat{\bm{g}}_{n} - \mathring{\bm{g}})) = A^{(m')} B_{n}(\hat{\bm{g}}_{n}^{*} - \mathring{\bm{g}})] = 1.
\end{eqnarray}

Next, since $\varphi$ is a continuous map, which follows from statement $(iii)$ of Lemma \ref{propisotinf}, and $A^{(m')}$ is a linear operator, then from the continuous mapping theorem it follows that 
\begin{eqnarray}\label{convphiABng}
A^{(m')} \varphi(B_{n}(\hat{\bm{g}}_{n} - \mathring{\bm{g}}))  \stackrel{d}{\to} A^{(m')}\varphi(\bm{\lambda}).
\end{eqnarray} 
and, using Lemma \ref{propconv} and result (\ref{eqaulAphiinf}), we prove
\begin{eqnarray*}\label{}
A^{(m')} B_{n}(\hat{\bm{g}}_{n}^{*} - \mathring{\bm{g}}) \stackrel{d}{\to} A^{(m')}\varphi(\bm{\lambda}).
\end{eqnarray*}

Note, that the number $m'$ is an arbitrary finite integer. Also, since $\varphi$ is a continuous map, then the law of $\varphi(\bm{\lambda})$ has the same continuity sets as $\bm{\lambda}$. Using Lemma \ref{convl2} we finish the proof of the theorem. \eop

Recall that the cardinality of any comparable level set $\mathcal{X}^{(v, l)}$ is less than infinity whenever $\mathring{g}_{v,l} \neq 0$. Then, as in the finite case, we note that the order constraints on $\mathcal{X}^{(v, l)}$ can be expressed in the form $\bm{A}\bm{g} \geq \bm{0}$, for some matrix $\bm{A}$. Therefore, one can use the results in \cite{Andrews} to describe the behaviour of $[\varphi(\bm{g})]^{(v, l)}$ when $|\mathcal{X}^{(v, l)}| < \infty$. It follows from Theorem 5 in \cite{Andrews} that the distribution of $[\varphi(\bm{g})]^{(v, l)}$ is a mixture of $2^{|\mathcal{X}^{(v, l)}|}$ distributions of the projections of $\bm{g}$ onto the cone $\bm{A}_{t}\bm{g} \geq \bm{0}$, where the matrices $\bm{A}_{t}$, for $t = 1, \dots, 2^{|\mathcal{X}^{(v, l)}|}$, are comprised of the rows of the matrix $\bm{A}$.

Next, let us consider the case of non-constant weights $\bm{w}$. In this section until now we assumed that the vector of weights satisfies the condition in (\ref{conditw}), it is  fixed, $\bm{w}_{n} = \bm{w}$, so it does not depend on $n$, and the random elements $\hat{\bm{g}}_{n}$ in Assumption \ref{asthetainf} all take their values in $(\bm{l}_{2}^{\bm{w}}, \mathcal{B})$, for some fixed $\bm{w}$, with $\mathcal{B}$ the Borel $\sigma$-algebra generated by the topology which is generated by the natural norm of $\bm{l}_{2}^{\bm{w}}$. 

Now we consider some non-random sequence $\{\bm{w}_{n}\}_{n\geq1}$, taking values in the space $\mathbb{R}^{\infty}$, where each $\bm{w}_{n}$ satisfies the condition in (\ref{conditw}). The sequence $\{\bm{w}_{n}\}_{n\geq1}$ converges in some norm $\norm{\cdot}_{R}$ on $\mathbb{R}^{\infty}$ to some non-random vector $\bm{w}$, which also satisfies the condition in (\ref{conditw}). Next, let $\mathcal{B}_{n}$ denotes the Borel $\sigma$-algebra generated by the topology which is generated by the natural norm in $\bm{l}_{2}^{\bm{w_{n}}}$. The next lemma shows that the normed spaces $\bm{l}_{2}^{\bm{w_{n}}}$ are all equivalent. 

\blem\label{equivlw}
Suppose that two vectors $\bm{w}_{1}$ and $\bm{w}_{2}$ satisfy the condition in (\ref{conditw}). Then the normed spaces $\bm{l}_{2}^{\bm{w_{1}}}$ and $\bm{l}_{2}^{\bm{w_{2}}}$ are equivalent. 
\elem

Therefore, since the normed spaces $\bm{l}_{2}^{\bm{w_{n}}}$ are all equivalent, then the topologies generated by these norms are the same. Then, the Borel $\sigma$-algebras $\mathcal{B}_{n}$ generated by these topologies are also the same. Therefore, the measurable spaces $(\bm{l}_{2}^{\bm{w}_{n}}, \mathcal{B}_{n})$ are all the same and we will suppress the index $n$. 

Next, analogously to the finite case, let us introduce the notation $\varphi^{\bm{w}}(\bm{\theta})$
\begin{eqnarray}\label{phiopwinf}
\varphi^{\bm{w}}(\bm{\theta})_{i} = \theta^{*\bm{w}(v',l')}(x_{i}),
\end{eqnarray} 
where $\theta^{*\bm{w}(v',l')}(x)$ is the isotonic regression with weights $\bm{w}$ of the restriction of $\theta(x)$ to the comparable level set $\mathcal{X}^{(v', l')}$ of $ \mathring{g}(x)$, where the indices $v'$ and $l'$ are such that $x_{i} \in \mathcal{X}^{(v', l')}$. The next theorem gives the limiting distribution of $\hat{\bm{g}}_{n}^{*\bm{w_{n}}}$.

\bth\label{thmasyminfw}
Suppose the Assumption \ref{asthetainf} holds. Then the asymptotic distribution of the isotonic regression $\hat{\bm{g}}_{n}^{*\bm{w_{n}}}$ of the basic estimator $\hat{\bm{g}}_{n}$ is given by
\begin{eqnarray}\label{asisot}
B_{n}(\hat{\bm{g}}_{n}^{*\bm{w}_{n}} - \mathring{\bm{g}}) \stackrel{d}{\to} \varphi^{\bm{w}}(\bm{\lambda}),
\end{eqnarray}
where $\varphi$ is the operator, defined in (\ref{phiopwinf}). \eth
\prf
First, we note that the result of Lemma \ref{convl2} holds, if we assume that the random vectors $\bm{Z}_{n}$, for $n = 1, \dots, \infty$ take their values in $\bm{l}_{2}^{\bm{w}_{n}}$, if all elements of $\bm{w}_{n}$ and its limit $\bm{w}$ satisfy the condition in (\ref{conditw}): This follows from the fact that the measurable spaces $(\bm{l}_{2}^{\bm{w}_{n}}, \mathcal{B}_{n})$ are equivalent, which was proved in Lemma \ref{equivlw}.

The rest of the proof is exactly the same as for Theorem \ref{thmasyminf} with $\varphi$  and $\hat{\bm{g}}_{n}^{*}$ suitable changed to $\varphi^{\bm{w}}$ and $\hat{\bm{g}}_{n}^{*\bm{w}_{n}}$. Also, recall that the result of Lemma \ref{isonlevelsinf} does not depend on the weights $\bm{w}_{n}$.
\eop

\section{Application to bimonotone probability mass function and regression function estimation, and extensions to $d$-dimensional problems}\label{sec:applications}
In this section we consider the problems of estimation of a bimonotone regression function, \tcb{in subsection \ref{subsec:bimon-reg}}{, and of a bimonotone probability mass function, \tcb{in subsection \ref{subsec:bimon-pmf}}. Also, we consider the generalisation to the case of $d$-dimensional support, \tcb{in subsection \ref{subsec:d-mon}}.

First, let us introduce a bimonotone order relation $\preceq$ on a set $\mathcal{X} := \{ \bm{x}= (i_{1} ,i_{2})^{T}: i_{1}=1,2,\dots, r_{1},  i_{2}=1,2,\dots, r_{2}\}$, with $r_{1},r_{2} \leq \infty$ in the following way. For any $\bm{x}_{1}$ and $\bm{x}_{2}$ in $\mathcal{X}$ we have $\bm{x}_{1} \preceq \bm{x}_{2}$ if and only if $x_{1,1} \leq x_{2,1}$ and $x_{1,2} \leq x_{2,2}$. The order relation $\preceq$ is a partial order, because it is reflexive, transitive, antisymmetric, but there are elements in $\mathcal{X}$ which are noncomparable. 

Second, note that $\mathcal{X}$ with the order relation $\preceq$ defined above is non-de\-com\-posable, because for any $\bm{x}_{1}=(x_{1,1}, x_{1,2})$ and $\bm{x}_{2}=(x_{2,1}, x_{2,2})$ in $\mathcal{X}$ there exist $\bm{x}_{3}=(x_{3,1}, x_{3,2}) \in \mathcal{X}$ such that $x_{3,1} \geq x_{1,1}$, $x_{3,1} \geq x_{2,1}$ and $x_{3,2} \geq x_{1,2}$, $x_{3,2} \geq x_{2,2}$, which means that  $\bm{x}_{1} \sim  \bm{x}_{3}$ and $\bm{x}_{2} \sim  \bm{x}_{3}$, or $\bm{x}_{4}=(x_{4,1}, x_{4,2}) \in \mathcal{X}$ such that $x_{4,1} \leq x_{1,1}$, $x_{4,1} \leq x_{2,1}$ and $x_{4,2} \geq x_{1,2}$, $x_{4,2} \geq x_{2,2}$, which means that  $\bm{x}_{1} \sim  \bm{x}_{4}$ and $\bm{x}_{2} \sim  \bm{x}_{4}$. Therefore, in a partition (\ref{partX}) $k=1$. Also, following our notations above we denote by $\mathcal{I}$ the set of indices of the domain $\mathcal{X}$ and use the same notation $\preceq$ for the order relation on $\mathcal{I}$ generated by $\mathcal{X}$.

A real valued function $g(\bm{x})$ is bimonotone increasing, i.e. isotonic with respect to bimonotone order relation $\preceq$ on a set $\mathcal{X}$, if whenever $\bm{x}_{1} \preceq \bm{x}_{2}$ one has $g(\bm{x}_{1}) \leq g(\bm{x}_{2})$, cf. \cite{dumbgen}. A real valued function $h(\bm{x})$ is called a bimonotone decreasing function, if whenever $\bm{x}_{1} \preceq \bm{x}_{2}$ one has $h(\bm{x}_{1}) \geq h(\bm{x}_{2})$. In the last case the function $h(\bm{x})$ is called antitonic with respect to the order relation $\preceq$ on a set $\mathcal{X}$, cf. \cite{barlowstatistical, robertsonorder}. Note that a fucntion $h(\bm{x})$ is antitonic if and only if $g(\bm{x}) = -h(\bm{x})$ is isotonic with respect to the order relation $\preceq$ on the set $\mathcal{X}$.

\subsection{Estimation of a bimonotone increasing regression function}\label{subsec:bimon-reg}
The problem of estimation of a bimonotone regression function via least squares was studied in detail in \cite{dumbgen}, where the authors described an algorithm for minimization of a smooth function under bimonotone order constraints.

Suppose we have observed $Z_{i} = (\bm{x}_{i}, Y_{i})$, $i = 1, \dots, n$, with $\bm{x}_{i}$ the design points taking values from the set $\mathcal{X} := \{ \bm{x}= (i_{1} , i_{2})^{T}: i_{1} = 1,2,\dots, r_{1},  i_{2} = 1,2,\dots, r_{2},\}$, with $r_{1}, r_{2} < \infty$ and $Y_{i}$ real valued random variables defined in the regression model
\begin{eqnarray*}
Y_{i} = \mathring{g}(\bm{x}_{i}) + \varepsilon_{i}, \quad i=1, \dots, n,
\end{eqnarray*}
where $\varepsilon_{i}$ is a sequence of identically distributed random variables with $\mathbb{E}[\varepsilon_{i}] = 0$, $\mathrm{Var}[\varepsilon_{i}] = \sigma^{2} < \infty$. 

Next, \cite{dumbgen} showed that the least squares estimator of $\mathring{g}(\bm{x})$ under bimonotone constraints is given by
\begin{eqnarray}\label{rbmr}
g_{n}^{*} = \underset{f \in \mathcal{F}^{is}}{\argmin} \sum_{\bm{x} \in \mathcal{X}}(f(\bm{x}) - \hat{g}_{n}(\bm{x}))^{2}w_{\bm{x}}^{(n)},
\end{eqnarray}
where $\mathcal{F}^{is}$ denotes the set of all bounded bimonotone increasing functions on $\mathcal{X}$, $\hat{g}_{n}(\bm{x})$ is the average of $Y_{i}$, $i =1, \dots, n$, over the design element $\bm{x}$, i.e. 
\begin{eqnarray}\label{unrbmr}
\hat{g}_{n}(\bm{x}) =  \frac{\sum_{i=1}^{n} Y_{i} 1 \{ \bm{x}_{i} = \bm{x} \}} {\sum_{i=1}^{n} 1 \{ \bm{x}_{i} = \bm{x} \}}
\end{eqnarray}
and 
\begin{eqnarray}\label{wxn}
w_{\bm{x}}^{(n)} = \frac{\sum_{i=1}^{n} 1\{\bm{x}_{i} = \bm{x}\}}{n}.
\end{eqnarray}

Note that $g_{n}(\bm{x})$ in (\ref{unrbmr}) is the unconstrained least squares estimate of $\mathring{g}(\bm{x})$. The asymptotic properties of nonlinear least squares estimators were studied in \cite{jennr, chienwu}. Assume that the design points $\bm{x}_{i}$, with $i = 1, \dots, n$, satisfy the following condition
\begin{eqnarray}\label{wxnconv}
\bm{w}^{(n)} \to \bm{w},
\end{eqnarray}
as $n \to \infty$, where $\bm{w}^{(n)}$ is a sequence of vectors in $\mathbb{R}_{+}^{r_{1} \times r_{2}}$ whose components are from (\ref{wxn}), and $\bm{w} \in \mathbb{R}_{+}^{r_{1} \times r_{2}}$. Given the condition in (\ref{wxnconv}) is satisfied, the basic estimator $\hat{g}_{n}(\bm{x})$ is consistent and has the following asymptotic distribution 
\begin{eqnarray}\label{asunreg}
n^{1/2}(\hat{\bm{g}}_{n} - \bm{g}) \stackrel{d}{\to} \bm{Y}_{\bm{0}, \Sigma},
\end{eqnarray}
where $\bm{Y}_{\bm{0}, \Sigma}$ is a Gausian vector with mean zero and diagonal covariance matrix $\Sigma$, whose elements are given by $\Sigma _{ii} = \sigma^{2} w_{i}$, for $i = 1, \dots, r \times s$, cf. Theorem 5 in \cite{chienwu}.  

We next derive the asymptotic distribution of the regression function under bimonotone constraints.
\bth\label{bimonotreg}
Given that the condition (\ref{wxnconv}) on the design points is satisfied, the asymptotic distribution of the regression function $\hat{\bm{g}}_{n}^{*}(\bm{x})$ under bimonotone constraints is given by
\begin{eqnarray}\label{asisot}
n^{1/2}(\hat{\bm{g}}_{n}^{*} - \mathring{\bm{g}}) \stackrel{d}{\to} \varphi^{\bm{w}}(\bm{Y}_{\bm{0}, \Sigma}),
\end{eqnarray}
where $\varphi^{\bm{w}}$ is the operator defined in (\ref{phiopw}) and $\bm{Y}_{\bm{0}, \Sigma}$ is a Gaussian vector defined in (\ref{asunreg}). \eth
\prf
The requirements of Assumption \ref{astheta} are satisfied. Therefore the result follows from  Theorem \ref{thmasyminf}.
\eop

Note, that the limit distribution $\varphi^{\bm{w}}(\bm{Y}_{\bm{0}, \Sigma})$ in (\ref{asisot}) is the concatenation of the separate weighted isotonic regressions of the restrictions of the Gaussian vector $\bm{Y}_{\bm{0}, \Sigma}$ to the comparable level sets of the regression function $\mathring{g}(\bm{x})$. 

\subsection{Estimation of a bimonotone decreasing probability mass function}\label{subsec:bimon-pmf}
In this subsection we treat the problem of estimating a bimonotone {\em decreasing } probability mass function on $\mathbb{Z}^{+}_{2}$. Note that this is a natural order restriction on the pmfs defined on $\mathbb{Z}^{+}_{2}$, since a positive bimonotone increasing function on $\mathbb{Z}^{+}_{2}$ does not belong to $l_{2}$. 

Suppose that we have observed $Z_{1}, Z_{2}, \dots, Z_{n}$ i.i.d. random variables taking values in $\mathcal{X} =\mathbb{Z}^{+}_{2}:= \{ (i_{1} ,i_{2})^{T}: i_{1} = 1,2,\dots, \infty, i_{2} = 1,2,\dots, \infty,\}$ with probability mass function $\bm{p}$. The empirical estimator of $\bm{p}$ is then given by
\begin{eqnarray}\label{unrMLEp}
	\hat{p}_{n, \bm{i}} = \frac{n_{\bm{i}}}{n}, \quad n_{\bm{i}}= \sum_{j=1}^{n} 1 \{ Z_{j} = \bm{x}_{\bm{i}} \}, \quad \text{$ \bm{i} \in \mathcal{I}$},
\end{eqnarray}
and it is also the unrestricted mle, which generally does not satisfy the bimonotonicity constraints introduced above. However, $\hat{\bm{p}}_{n}$ is consistent, i.e. $\hat{\bm{p}}_{n} \stackrel{p}{\to} \bm{p}$ and asymptotically Gaussian
\begin{eqnarray}\label{empasym}
n^{1/2}(\hat{\bm{p}}_{n} - \bm{p}) \stackrel{d}{\to} \bm{Y}_{\bm{0}, C},
\end{eqnarray}
where $\bm{Y}_{\bm{0}, C}$ is a Gaussian process in $l_{2}$, with mean zero and the covariance operator $C$ such that $\langle C\bm{e}_{\bm{i}},\bm{e}_{\bm{i'}} \rangle  = p_{\bm{i}} \delta_{\bm{i},\bm{i'}} - p_{\bm{i}}p_{\bm{i'}}$, with $\bm{e}_{\bm{i}} \in \bm{l}_{2}$ the orthonormal basis in $\bm{l}_{2}$ such that in a vector $\bm{e}_{\bm{i}} $ all elements are equal to zero but the one with the index $\bm{i}$ is equal to 1, and $\delta_{ij}=1$, if $i=j$ and $0$ otherwise, cf. \cite{jankowski2009estimation}.

The constrained mle $\hat{\bm{p}}^{*}_{n}$ of $\bm{p}$ is then given by the isotonic regression of the empirical estimator $\hat{\bm{p}}_{n}$ over the set $\mathcal{X}$
with respect to the pre-order $\preceq$
\begin{eqnarray}\label{bimpmfisot}
\hat{\bm{p}}^{*}_{n} = \underset{\xi \in \mathcal{F}^{an}}{\argmin} \sum_{\bm{x} \in \mathcal{X}}(\xi_{\bm{x}} - \hat{p}_{n,\bm{x}})^{2},
\end{eqnarray}
where $\mathcal{F}^{an}$ denotes the set of all bimonotone decreasing (antitonic with respect to $\preceq$) functions on $\mathcal{X}$. This result is shown on pages 45--46 in \cite{barlowstatistical} and pages 38--39 in \cite{robertsonorder}. 

Next we make the following substitution
\begin{eqnarray}\label{thetmin}
\bm{\theta} &=& - \bm{p},\nonumber\\
\hat{\bm{\theta}}_{n} &=& - \hat{\bm{p}}_{n},\\\nonumber
\hat{\bm{\theta}}^{*}_{n} &=& - \hat{\bm{p}}^{*}_{n}. 
\end{eqnarray}
Therefore $\hat{\bm{\theta}}^{*}_{n}$ is the isotonic regression of $\hat{\bm{\theta}}_{n}$, i.e.
\begin{eqnarray}\label{bimpmfisot1}
\hat{\bm{\theta}}^{*}_{n} = \underset{\xi \in \mathcal{F}^{is}}{\argmin} \sum_{\bm{x} \in \mathcal{X}}(\xi_{\bm{x}} - \hat{\theta}_{n,\bm{x}})^{2},
\end{eqnarray}
where $\mathcal{F}^{is}$ denotes the set of all bimonotone increasing (isotonic with respect to $\preceq$) functions on $\mathcal{X}$.

We next derive the asymptotic distribution of the bimonotone mle $\hat{\bm{p}}^{*}_{n}$ as a corollary of Theorem \ref{thmasyminf}.
\bth\label{bimonot}
The asymptotic distribution of the constrained mle $\hat{\bm{p}}^{*}_{n}$ of a bimonotone probability mass function $\bm{p}$ is given by
\begin{eqnarray}\label{asisotpmf}
n^{1/2}(\hat{\bm{p}}_{n}^{*} - \bm{p}) \stackrel{d}{\to} \varphi(\bm{Y}_{\bm{0}, C}),
\end{eqnarray}
where $\varphi$ is the operator defined in (\ref{phiinf}) and $\bm{Y}_{\bm{0}, C}$ is a Gaussian process in $l_{2}$ defined in (\ref{empasym}). \eth
\prf
The requirements of Assumption \ref{asthetainf} for the sequence $\hat{\bm{\theta}}_{n}$ defined in (\ref{thetmin}) are satisfied. Therefore from  Theorem \ref{thmasyminf} it follows that
\begin{eqnarray*}\label{}
n^{1/2}(\hat{\bm{\theta}}_{n}^{*} - \bm{\theta}) \stackrel{d}{\to} \varphi(\bm{Y}_{\bm{0}, C})
\end{eqnarray*}
and using the substitution (\ref{thetmin}) we finish the proof.
\eop

Note, that the limit distribution $\varphi(\bm{Y}_{\bm{0}, C})$ in (\ref{asisotpmf}) is the concatenation of the separate isotonic regressions of the restrictions of the Gaussian vector $\bm{Y}_{\bm{0}, C}$, defined in (\ref{empasym}), to the comparable level sets of the true pmf $\bm{p}$. In the case of the estimation of a one-dimensional decreasing pmf we get the result previously obtained in Theorem 3.8 in \cite{jankowski2009estimation}, i.e. the limit distribution $\varphi(\bm{Y}_{\bm{0}, C})$ in (\ref{asisotpmf}) is the concatenation of the separate isotonic regressions of the restrictions of the Gaussian vector $\bm{Y}_{\bm{0}, C}$ to the constant regions of the true pmf $\bm{p}$.

\subsection{Generalisation to the case of $d$-dimensional monotone functions}\label{subsec:d-mon}
The results obtained in Theorems \ref{bimonotreg} and \ref{bimonot} can be directly generalised to the case of estimation of a $d$-dimensional monotone ($d$-monotone) regression function and a $d$-monotone pmf.  

Let us consider a set 
\begin{eqnarray}\label{Xdmon}
\mathcal{X} := \{ \bm{x}= (i_{1} ,i_{2}, \dots, i_{d})^{T}: i_{1}=1,2,\dots, r_{1},  i_{2}=1,2,\dots, r_{2}, \dots, \nonumber\\ i_{d}=1,2,\dots, r_{d}\}, \text{ with } d < \infty, r_{1},r_{2}, \dots, r_{d} \leq \infty
\end{eqnarray}
and introduce a $d$-monotone order relation $\preceq$ on it in the following way. For any $\bm{x}_{1}$ and $\bm{x}_{2}$ in $\mathcal{X}$ we have $\bm{x}_{1} \preceq \bm{x}_{2}$ if and only if $x_{1,1} \leq x_{2,1}, x_{1,2} \leq x_{2,2}, \dots, x_{1,d} \leq x_{2,d}$. Similarly to the bimonotone case, it can be shown that the order relation $\preceq$ is a partial order and $\mathcal{X}$ is non-decomposable.

Suppose we have observed $Z_{i} = (\bm{x}_{i}, Y_{i})$, $i = 1, \dots, n$, with $\bm{x}_{i}$ the design points taking values from the set $\mathcal{X}$ defined in (\ref{Xdmon}), with $r_{1}, r_{2}, \dots,  r_{d} < \infty$ and $Y_{i}$ real valued random variables defined in the regression model
\begin{eqnarray*}
Y_{i} = \mathring{g}(\bm{x}_{i}) + \varepsilon_{i}, \quad i=1, \dots, n,
\end{eqnarray*}
where $\varepsilon_{i}$ is a sequence of identically distributed random variables with $\mathbb{E}[\varepsilon_{i}] = 0$, $\mathrm{Var}[\varepsilon_{i}] = \sigma^{2} < \infty$. 

The least squares estimate of $\mathring{g}(\bm{x})$ under $d$-monotone constraints is given by
\begin{eqnarray*}\label{rdmr}
g_{n}^{*} = \underset{f \in \mathcal{F}^{is}}{\argmin} \sum_{\bm{x} \in \mathcal{X}}(f(\bm{x}) - \hat{g}_{n}(\bm{x}))^{2}w_{\bm{x}}^{(n)},
\end{eqnarray*}
where $\mathcal{F}^{is}$ denotes the set of all bounded $d$-monotone functions on $\mathcal{X}$, the expressions for $\hat{g}_{n}(\bm{x})$ and $w_{\bm{x}}^{(n)}$ are the same as in bimonotone case, i.e. given in (\ref{unrbmr}) and (\ref{wxn}), respectively. Therefore, under condition (\ref{wxnconv}) on the design points $\bm{x}_{i}$, we obtain the following corollary.

\bcor\label{dmonotreg}
The asymptotic distribution of the regression function $\hat{\bm{g}}_{n}^{*}(\bm{x})$ under $d$-monotone constraints is given by
\begin{eqnarray*}\label{}
n^{1/2}(\hat{\bm{g}}_{n}^{*} - \mathring{\bm{g}}) \stackrel{d}{\to} \varphi^{\bm{w}}(\bm{Y}_{\bm{0}, \Sigma}),
\end{eqnarray*}
where $\varphi^{\bm{w}}$ is the operator defined in (\ref{phiopw}) and $\bm{Y}_{\bm{0}, \Sigma}$ is a Gaussian vector defined in (\ref{asunreg}). \ecor
\prf
The requirements of Assumption \ref{astheta} are satisfied. Therefore the result follows from  Theorem \ref{thmasyminf}.
\eop

Next suppose that we have observed $Z_{1}, Z_{2}, \dots, Z_{n}$ i.i.d. random variables taking values in $\mathcal{X}$ defined in (\ref{Xdmon}), with $r_{1}, r_{2}, \dots,  r_{d} \leq \infty$ with probability mass function $\bm{p}$. The mle $\hat{\bm{p}}^{*}_{n}$ of $\bm{p}$ with $d$-monotone decreasing constraints is then given by 
\begin{eqnarray*}\label{dmpmfisot}
\hat{\bm{p}}^{*}_{n} = \underset{\xi \in \mathcal{F}^{an}}{\argmin} \sum_{\bm{x} \in \mathcal{X}}(\xi_{\bm{x}} - \hat{p}_{n,\bm{x}})^{2},
\end{eqnarray*}
where $\hat{\bm{p}}_{n}$ is the empirical estimator defined in (\ref{unrMLEp}), $\mathcal{F}^{an}$ denotes the set of all $d$-monotone decreasing functions on $\mathcal{X}$. The asymptotic distribution of $\hat{\bm{p}}^{*}_{n}$ is given in the following corollary.
\bcor\label{dmonot}
The asymptotic distribution of the constrained mle $\hat{\bm{p}}^{*}_{n}$ of a $d$-monotone probability mass function $\bm{p}$ is given by
\begin{eqnarray}\label{asisot}
n^{1/2}(\hat{\bm{p}}_{n}^{*} - \bm{p}) \stackrel{d}{\to} \varphi(\bm{Y}_{\bm{0}, C}),
\end{eqnarray}
where $\varphi$ is the operator defined in (\ref{phiinf}) and $\bm{Y}_{\bm{0}, C}$ is a Gaussian process in $l_{2}$ defined in (\ref{empasym}). \ecor
\prf
Making the same substitution as in a bimonotone case, i.e. in (\ref{thetmin}) we note that the requirements of Assumption \ref{asthetainf} are satisfied. Therefore the result follows from Theorem \ref{thmasyminf}.
\eop

\section{Conclusions and discussion}\label{sec:discussion}
We have derived the limit distribution of an estimator that is obtained as the $l^2$-projection of a basic preliminary estimator on the space of functions that are defined on a countable set, and that are monotone with respect to a pre-order on that countable set. Immediate applications that we have stated results for are to the estimation of $d-$monotone pmfs and regression functions.

We would like to emphasize a qualitative difference between the estimation of a pmf over a subset of ${\mathbb Z}_+^d$ and the estimation of a pdf over a subset of ${\mathbb R}^d_+$. We note first that limit distribution results for monotone pdf estimators, to our knowledge, exist only for the case $d=1$, cf. however \cite{polonik:1998} for the limit distribution of the non-parametric maximum likelihood estimator (npmle) of a bimonotone pdf (so when $d=2$), indexed by (the Lebesgue measure of) lower layers. For the case $d=1$, the isotonic regression estimator of a pdf is, for independent data, equivalent to the npmle, i.e. the Grenander estimator, and for dependent data does not have the interpretation of an npmle, cf. \cite{anevhos} for the limit distribution results for the monotone restricted pdf estimator for arbitrary dependence assumptions on the data. The limit distribution in the independent data case is then the well known Chernoff distribution mentioned above, and for dependent data different, cf.  Theorem 10 $(ii)$ and Theorem 11 in \cite{anevhos}.

Note also that, in the case $d=1$, the order restricted estimator of a pdf is a {\em local} estimator, in the sense that it uses data in a shrinking neighbourhood around the point of interest, say $t_0\in{\mathbb R}$, to calculate the value of the pdf at $t_0$, and the size of the neighbourhood is of the order $n^{-1/3}$ for independent data, and of a different order for dependent data, cf. Table 1 of Section 5 in \cite{anevhos} for an overview of the possible orders related to the dependence of the data. Any sensible estimator of a monotone pdf for $d\geq 2$ will also use data in a shrinking neighbourhood around the point of interest, cf. e.g. \cite{han:wang:chatterjee:samworth:2017} for a discussion about rates in this connection. Furthermore, as argued e.g. in  \cite{han:wang:chatterjee:samworth:2017}, the rates in higher dimensions are slower for monotone pdf estimation. This is in sharp contrast to the problems treated in this paper, on monotone pmf estimation, and is explained by the fact that the resulting estimator for those problems is a {\em global} estimator, i.e. it uses data points in a set of size $O(1)$ around the point of interest to obtain the estimator, irrespective of the dimension $d$. The fact that estimators of pdf are local and of pmf are global, also accounts for that one is able to obtain process limit distribution results for the pmf estimator, whereas it is only possible to obtain pointwise limit distribution results for the pdf estimators. Similar qualitative differences exist between the continuous case regression problem, treated e.g. in \cite{anevhos} for $d=1$, and the discrete regression problem treated in this paper, see also \cite{dumbgen} for the discrete bimonote regression problem.

The results stated in this paper are general in terms of the demands on the basic estimator and on the underlying empirical process. In fact, Assumptions \ref{astheta} and \ref{asthetainf} only require that there is a limit process for the basic estimator, and do not specify any requirements of e.g. dependence for the data. Note, however, that if one does not require independence of the data, then the identity between the isotonic regression of a pmf and the mle of a pmf vanishes, since the product of the marginal pmfs is then not the full likelihood.

By allowing dependent data, we are in a position to straight-forwardly obtain limit distributions  in general situations. One problem that comes to mind is that of isotonic regression of an ordered pmf on a DAG. The assumption of monotonicity of the pmf with respect to the natural tree order on the DAG is sensible; one can e.g. imagine the DAG describing the, say, three categories that may influence the monthly salary of an employee at a large facility, with the DAG structure given by the  (matrix) pre-order on the three categories. Then, given data on employees salary and covariate readings for the three categories, one may first construct the empirical estimate of the pmf and next isotonize that. Knowing the limit distribution of the empirical estimator immediately gives us the limit distribution of the isotonized estimator, irrespective of whether data are independent or not.

Finally we note that the limit distribution of an isotonized estimator depends on the structure of the true parameter, which is unknown in general. Therefore the limit distribution result's use in practice may be limited, and this is true also for the results obtained in \cite{jankowski2009estimation}. In a next paper we use the results of this paper to perform model selection in the estimation of a decreasing pmf. Note however, that analogously to the unimodal pmf estimation problem in \cite{baljan}, the error reduction property of the isotonic regression can be used for the construction of conservative confidence bands for pmf or monotone regression function. 

\section{Appendix}
\textbf{Proof of Lemma \ref{propisot}}.
The statements $(i)$, $(ii)$, $(iii)$ and $(iv)$ are from \cite{robertsonorder} (Theorems 1.3.1, 1.3.3, 1.4.4 and 1.3.4). The statements $(v)$ and $(vi)$ are proved in \cite{barlowstatistical} (Theorems 2.2 and 1.8).

Note that statement $(ii)$ means that if the basic estimator $\hat{\bm{g}}_{n}$ satisfies a linear restriction, e.g. $\sum_{i=1}^{s} w_{i} \hat{g}_{n, i} = c$, with some positive reals $w_{i}$, then the same holds for its isotonic regression with the weights $\bm{w}$, i.e. for $\hat{\bm{g}}^{*}_{n}$ one has $\sum_{i=1}^{s} w_{i} \hat{g}^{*}_{n, i} = c$.  
\eop

\textbf{Proof of Lemma \ref{isotoutpartlem}}.
Let $g(x)$ be an arbitrary real-valued function defined on $ \mathcal{X}$. From the definition of the isotonic regression 
\begin{eqnarray*}
g^{*}&=&\underset{f \in \mathcal{F}^{is}}{\argmin} \sum_{x \in \mathcal{X}}(f(x) - g(x))^{2}w_{x}\\
&=&\underset{f \in \mathcal{F}^{is}}{\argmin}\sum_{v=1}^{k} \sum_{x \in \mathcal{X}^{(v)}}(f(x) - g(x))^{2}w_{x}\\
&=&\sum_{v=1}^{k}\underset{f^{(v)} \in \mathcal{F}^{is}_{(v)}}{\argmin}\sum_{x \in \mathcal{X}^{(v)}}(f^{(v)}(x) - g^{(v)}(x))^{2}w_{x},
\end{eqnarray*}
where $f^{(v)}$ is the restriction of the function $f: \mathcal{X} \to \mathbb{R}$ to the set $\mathcal{X}^{(v)}$. The second equality follows from ($\ref{partX}$) and the last equality follows from the fact that since the elements from the different partition sets $\mathcal{X}^{(v)}$ are noncomparable, then any function $f \in \mathcal{F}^{is}$ can be written as a concatenation of $f^{(v)} \in \mathcal{F}^{is}_{(v)}$, with no restrictions imposed on the values of $f^{(v_{1})}$ and $f^{(v_{2})}$ for $v_{1} \neq v_{2}$.
\eop

\textbf{Proof of Lemma \ref{isonlevels}}.
First, note that if the condition of the lemma is satisfied, then the function $g^{*}(x)$ defined in (\ref{isregndcs}) on the set $\mathcal{X}$ is isotonic. This follows from Lemma \ref{propisot}, statement $(iv)$. Second, assume that the function $g^{*}(x)$ defined in (\ref{isregndcs}) is not an isotonic regression of $g(x)$. This means that there exists another function $\tilde{g}(x)$, such that
\begin{eqnarray}\label{ineqprlem}
\sum_{x \in \mathcal{X}}(\tilde{g}(x) - g(x))^{2}w_{x} < \sum_{x \in \mathcal{X}}(g^{*}(x) - g(x))^{2}w_{x},
\end{eqnarray}
Using the partition of $\mathcal{X}$, (\ref{ineqprlem}) can be rewritten as 
\begin{eqnarray*}
\sum_{l=1}^{m} \sum_{x \in \mathcal{X}^{(l)}}(\tilde{g}(x) - g(x))^{2}w_{x} < \sum_{l=1}^{m} \sum_{x \in \mathcal{X}^{(l)}}(g^{*}(x) - g(x))^{2}w_{x}.
\end{eqnarray*}
Therefore, for some $l'$ we must have
\begin{eqnarray*}
\sum_{x \in \mathcal{X}^{(l')}}(\tilde{g}(x) - g(x))^{2}w_{x} < \sum_{x \in \mathcal{X}^{(l')}}(g^{*}(x) - g(x))^{2}w_{x}
\end{eqnarray*}
or, equivalently,
\begin{eqnarray*}
\sum_{x \in \mathcal{X}^{(l')}}(\tilde{g}^{(l')}(x) - g^{(l')}(x))^{2}w_{x} < \sum_{x \in \mathcal{X}^{(l')}}(g^{*(l')}(x) - g^{(l')}(x))^{2}w_{x},
\end{eqnarray*}
with $g^{(l')}(x)$, $\tilde{g}^{(l')}(x)$ and $g^{*(l')}(x)$ the restrictions to the comparable level set $\mathcal{X}^{(l')}$ of $g(x)$, $\tilde{g}(x)$ and $g^{*}(x)$, respectively. Since the function $g^{*(l')}(x)$ is the isotonic regression of the function $g^{(l')}(x)$ on the set $\mathcal{X}^{(l')}$, the last inequality contradicts the property of the uniqueness and existence of the isotonic regression $g^{*(l')}(x)$ (statement $(i)$ of Lemma \ref{propisot}).
\eop

\textbf{Proof of Lemma \ref{propconv}}.
This result follows from Theorem 3.1 in \cite{Billingsley} \eop

\textbf{Proof of Lemma \ref{propisotinf}}.
Statements $(i)$, $(ii)$ and $(iii)$ follow from Theorem 8.2.1, Corollary B of Theorem 8.2.7 and Theorem 8.2.5, respectively, in \cite{robertsonorder}, statements $(iv)$, $(v)$ and $(vi)$ follow from Corollary B of Theorem 7.9, Theorems 2.2 and Theorems 7.5 and 7.8, respectively, in \cite{barlowstatistical}.
\eop

\textbf{Proof of Lemma \ref{propisotinf}}.
The proof is exactly the same as in the finite case (Lemma \ref{isotoutpartlem}).

\eop

\textbf{Proof of Lemma \ref{isonlevelsinf}}.
The proof is exactly the same as in the case of a finite support (Lemma \ref{isonlevels}).
\eop

\textbf{Proof of Lemma \ref{convl2}}.
The space $\bm{l}_{2}^{\bm{w}}$ is separable and complete. Then, from Prokhorov's theorem \cite{Shrv2007}, it follows that the sequence $\bm{Z}_{n}$ is relatively compact, which means that every sequence from $\bm{Z}_{n}$ contains a subsequence, which converges weakly to some vector $\bm{Z}$. If the limits of the convergent subsequences are the same, then the result of the lemma holds.

Since the space $\bm{l}_{2}^{\bm{w}}$ is separable, the Borel $\sigma$-algebra equals the $\sigma$-algebra generated by open balls in $\bm{l}_{2}^{\bm{w}}$ \cite{Bgch2007}. Therefore, it is enough to show that the limit laws agree on the finite intersections of the open balls, since the finite intersections of open balls in $\bm{l}_{2}^{\bm{w}}$ constitute a $\pi$-system. 


Let us consider two arbitrary balls $B(\bm{z}_{1}, \varepsilon_{1})$ and $B(\bm{z}_{2}, \varepsilon_{2})$ in $\bm{l}_{2}^{\bm{w}}$ and note that 
\begin{eqnarray*}
B(\bm{z}_{1}, \varepsilon_{1}) &=&  \cap_{M_{1}\geq 1}C_{1,M_{1}},
\end{eqnarray*}
where $C_{1,M_{1}}$ is the following cylinder set in $\bm{l}_{2}^{\bm{w}}$
\begin{eqnarray*}
   C_{1,M_{1}}&=&\{\bm{y} \in \bm{l}_{2}^{\bm{w}}: \underset{j \in {\tilde{i}_{1}, \dots, \tilde{i}_{M_{1}}}}{\sum}|z_{1,j} - y_{j}|^{2}w_{j} < \varepsilon_{1}^{2}\}
\end{eqnarray*}
and
\begin{eqnarray*}
B(\bm{z}_{2}, \varepsilon_{2}) &=&  \cap_{M_{2}\geq 1}C_{2,M_{2}},
\end{eqnarray*}
where $C_{2,M_{2}}$ is the following cylinder set in $\bm{l}_{2}^{\bm{w}}$
\begin{eqnarray*}
   C_{2,M_{2}}&=&\{\bm{y} \in \bm{l}_{2}^{\bm{w}}: \underset{j \in {\tilde{i}_{1}, \dots, \tilde{i}_{M_{2}}}}{\sum}|z_{2,j} - y_{j}|^{2}w_{j} < \varepsilon_{2}^{2}\} 
\end{eqnarray*}
where the indices $\tilde{i}_{1}, \dots, \tilde{i}_{M}$ are the first $M$ indices from $\tilde{\mathcal{I}}$ and the index set $\tilde{\mathcal{I}}$ is the same for both $C_{1,M}$ and $C_{2,M}$.

Next, using the  associativity of the intersection, we obtain 
\begin{eqnarray}\label{ballsint}
B(\bm{z}_{1}, \varepsilon_{1}) \cap B(\bm{z}_{2}, \varepsilon_{2}) = \cap_{M\geq 1}(C_{1,M}\cap C_{2,M}).
\end{eqnarray}


The sequence of vectors $\bm{\tilde{Z}}^{(1, M)}_{n}$ converges weakly to $\bm{\tilde{Z}}^{(1, M)}$ for all finite M, therefore any subsequence of $\bm{\tilde{Z}}^{(1, M)}_{n}$ converges weakly to $\bm{\tilde{Z}}^{(1, M)}$. 
Therefore, by the continuity property of a probability measure, we have
\begin{eqnarray*}
      \mathbb{P}(B(\bm{z}_{1}, \varepsilon_{1}) \cap B(\bm{z}_{2}, \varepsilon_{2}))&=&\mathbb{P}( \cap_{M\geq 1}(C_{1,M}\cap C_{2,M}))\\
\lim_{M\to \infty}\mathbb{P}(C_{1,M}\cap C_{2,M})&=&\lim_{M\to \infty}\mathbb{P}^{(M)}(C_{1,M}\cap C_{2,M})
\end{eqnarray*}
where $\mathbb{P}^{(M)}$ denotes the law of $\bm{\tilde{Z}}^{(1, M)}$ and recall that  $\mathbb{P}^{(M)}$ is the same for all convergent subsequences of $\{\bm{Z}_{n}\}_{n\geq 1}$. Therefore, $\mathbb{P}(B(\bm{z}_{1}, \varepsilon_{1}) \cap B(\bm{z}_{2}, \varepsilon_{2}))$ are equal for all convergent subsequences. 

Thus, we have shown that the limit laws of the convergent subsequences of $\{ \bm{Z}_{n} \}_{n\geq 1}$ agree on the intersection of open balls, which proves the lemma. \eop

\textbf{Proof of Lemma \ref{equivlw}}.
First, we prove that if $\bm{w}$ satisfies the condition in (\ref{conditw}), then $\bm{x} \in \bm{l}_{2}^{\bm{w}}$ if and only if $\bm{x} \in \bm{l}_{2}$ ($\bm{l}_{2}$ is the space of all square summable sequences, i.e. $\bm{w} = (1, 1, \dots )$). Let $\bm{x} \in \bm{l}_{2}^{\bm{w}}$, then $\sum_{i=1}^{\infty} x_{i}^{2} w_{i} < \infty$ and we have
\begin{eqnarray*}\label{}
(\underset{i}{\inf} \{ w_{i} \}) \sum_{i=1}^{\infty} x_{i}^{2} \leq \sum_{i=1}^{\infty} x_{i}^{2} w_{i} < \infty.
\end{eqnarray*}

Therefore, since $\underset{i}{\inf} \{ w_{i} \} > 0$, we have that $\sum_{i=1}^{\infty} x_{i}^{2} < \infty$, which means that $\bm{x} \in \bm{l}_{2}$. 

Next, let $\bm{x} \in \bm{l}_{2}$, then $\sum_{i=1}^{\infty} x_{i}^{2} < \infty$ and we have
\begin{eqnarray*}\label{}
\sum_{i=1}^{\infty} x_{i}^{2} w_{i} \leq (\underset{i}{\sup} \{ w_{i} \}) \sum_{i=1}^{\infty} x_{i}^{2} < \infty,
\end{eqnarray*}
since $\underset{i}{\sup} \{ w_{i} \} < \infty$. Therefore, $\bm{x} \in \bm{l}_{2}^{\bm{w}}$.

Second, let $\norm{\cdot}_{\bm{w}}$ and $\norm{\cdot}$ denote the natural norms in $\bm{l}_{2}^{\bm{w}}$ and $\bm{l}_{2}$. We can prove that if $\bm{w}$ satisfies the condition in (\ref{conditw}), then $\bm{l}_{2}^{\bm{w}}$ and $\bm{l}_{2}$ are equivalent, i.e. there exist two positive constants $c_{1}$ and $c_{2}$ such that 
\begin{eqnarray}\label{equivnorms}
c_{1}\norm{\bm{x}} \leq \norm{\bm{x}}_{\bm{w}} \leq c_{2}\norm{\bm{x}},
\end{eqnarray}
if, for example, $c_{1} = \underset{i}{\inf} \{ w_{i} \}$ and $c_{2} = \underset{i}{\sup} \{ w_{i} \}$. Therefore, since the equivalence of norms is transitive, then  $\bm{l}_{2}^{\bm{w_{1}}}$ and $\bm{l}_{2}^{\bm{w_{2}}}$ are equivalent, provided $\bm{w_{1}}$ and $\bm{w_{2}}$ satisfy the condition in (\ref{conditw}).
\eop


\section*{Acknowledgements}
VP's research is fully supported and DA's research is partially supported by the Swedish Research Council, whose support is gratefully acknowledged.


\end{document}